\documentclass{amsart}
\usepackage[foot]{amsaddr}
\usepackage{graphicx,natbib,covington,enumerate,setspace,bussproofs,mathrsfs,hyperref,color} % Required for inserting images

\newtheorem{theorem}{Theorem}

\title[Artifical intelligence and inherent mathematical difficulty]{Artifical intelligence and \\ inherent mathematical difficulty}
\author{Walter Dean}
\address{Department of Philosophy, University of Warwick}
\email{W.H.Dean@warwick.ac.uk}
\author{Alberto Naibo}
\address{Department of Philosophy, Universit\'e Paris 1 Panth\'eon-Sorbonne}
\email{Alberto.Naibo@univ-paris1.fr}

\date{}

\begin{document}

\begin{abstract}  This paper explores the relationship of artificial intelligence to the task of resolving open questions in mathematics. We first present an updated version of a traditional argument that limitative results from computability and complexity theory show that proof discovery is an \textsl{inherently difficult problem}.  We then illustrate how several recent applications of artificial intelligence-inspired methods -- respectively involving automated theorem proving, \textsc{Sat}-solvers, and large language models  -- do indeed raise novel questions about the nature of mathematical proof.  We also argue that the  results obtained by such techniques do not tell against our basic argument.  This is so because they are embodiments of \textsl{brute force search} and are thus capable of deciding only statements of \textsl{low logical complexity}.
\end{abstract}

\maketitle

\vspace{2ex}

{\footnotesize
\begin{spacing}{1}
\noindent Suppose $\ldots$ that we could find a finite system of rules which enabled us to say whether any given formula was demonstrable or not. This system would embody a theorem of metamathematics. There is of course no such theorem and this is very fortunate, since if there were we should have a mechanical set of rules for the solution of all mathematical problems, and our activities as mathematicians would come to an end.  \hfill \citep[p. 16]{Hardy1929} 
\end{spacing}
}

\vspace{2ex}

{\footnotesize
\begin{spacing}{1}
\noindent On the basis of $\ldots$ the speed with which research in this field is progressing, I am willing to make the following predictions $\ldots$:
%\vspace{-.5ex}
\begin{enumerate}[1)]
\item That within ten years a digital computer will be the world's chess champion $\ldots$
%\vspace{-.7ex}
\item That within ten years a digital computer will discover and prove an important new mathematical theorem.  \hspace*{48ex}  \citep[p. 7]{Simon1958} 
\end{enumerate} 
\end{spacing}
}

\vspace{2ex}

\begin{footnotesize} 
\begin{spacing}{1}
\noindent [W]e are witnessing a steady increase in the intelligence of theorem proving software $\ldots$ The time will come when such crushers as Riemann's hypothesis and Goldbach's conjecture will be fair game for automated reasoning programs. For those of us who arrange to stick around, endless fun awaits us in the automated development and eventual enrichment of the corpus of mathematics.  
\hspace*{1ex} \hfill \citep[pp. 119-120]{Quaife1992}
\end{spacing}
\end{footnotesize}

\vspace{2ex}

\begin{footnotesize} 
\begin{spacing}{1}
\noindent A computer program $\ldots$ has come up with a major mathematical proof that would have been called creative if a human had thought of it. In doing so, the computer has, for the first time, got a toehold into pure mathematics $\ldots$ And the implications, some say, are profound, showing just how powerful computers can be at reasoning itself, at mimicking the flashes of logical insight or even genius that have characterized the best human minds. %\P \ Computers have found proofs of mathematical conjectures before, of course, but those conjectures were easy to prove. The difference this time is that the computer has solved a conjecture that stumped some of the best mathematicians for 60 years. 
\hfill \citep{Kolata1996}
\end{spacing}
\end{footnotesize}

\section{Introduction}
\label{sec1}

We write at a time when many endeavors are being called upon to reflect on the challenges brought about by recent developments in artificial intelligence.  It need hardly be pointed out that technologies falling under this umbrella have made considerable strides in a wide range of applications.   As a result, many  disciplines are in the midst of coming to terms with how artificial intelligence may come to affect how they are practiced or even their core subject matter.

This includes mathematics.    Mathematicians have, of course, employed computational methods since the beginning.   This includes the application of algorithms for  numerical, algebraic, and geometric calculation, some of which predate mechanical computers by thousands of years.  In the second half of the twentieth century, it came to include the use of software implementing numerical methods, computational algebra, and  visualization techniques.  Within the past twenty years, it has grown to subsume the use of automated systems to verify long or complex proofs originally obtained in outline by human mathematicians by formalizing them in first- or higher-order logic.   Within the past ten years initiatives have been proposed for making formalization in this style a routine part of mathematical practice. 

On the other hand, much of the current discussion of artificial intelligence has centered around techniques falling under the sub-umbrellas of \textsl{machine learning} and \textsl{generative artificial intelligence}.    At least  at the level of family resemblance, methods in this class are unlike those which  have previously been employed in mathematics in the sense that they embody statistical techniques for probabilistic prediction rather than logical techniques for deductive inference.   Nonetheless it is this class of technologies which are currently being discussed in regard to their ability to perform tasks which have traditionally been understood as the exclusive province of human insight or creativity.

The question thus arises whether  systems currently under investigation or ones which may evolve out of them will substantially transform the practice of mathematics as we now understand it.  In fact we are living at a moment when one might reasonably wonder whether the task of resolving open questions in mathematics by discovering appropriate proofs or refutations will succumb to routine automation.

This is a timely question with which we believe philosophers of mathematics ought to be engaging.\footnote{The papers collected in a recent volume of the \textsl{Bulletin of the American Mathematical Society} \citep{Adem2024} address the prospects of automation from the perspective of working mathematicians in a manner largely complementary to the argument which we advanced below.}    But it is also difficult to approach in a manner which avoids stipulating which methods exemplify ``artificial intelligence'' or prognosticating about the future course of technology.   These hazards notwithstanding,  the aim of this paper will be to argue that traditional limitative results from logic and theoretical computer science still have something to tell us about what we might hope to achieve.    In particular, we will argue that there are principled reasons to believe that resolving open questions in mathematics will remain difficult for computers for essentially the same reason it is difficult for humans -- i.e. because it is an inherently difficult task for any computational system to perform.

We will lay out our case as follows.   In \S 2, we will describe what we will mean by \textsl{proof discovery} and explain the sense in which we take it to be a central aspect of mathematical activity.  In \S 3, we will present what we will call the \textsl{basic argument} that proof discovery is inherently difficult. As we will see, a form of this argument was anticipated early in the history of artificial intelligence grounded in observations from computability and computational complexity theory.   Our task will thus be to illustrate how a refinement remains relevant today.      

In \S 4 we will examine three examples which might be taken to tell against either the correctness or significance of the basic argument.   These will take the form of instances in which techniques from artificial intelligence -- understood in the broad sense of \citep{Russell2021} -- have been successfully applied to resolve open mathematical questions.   The first involves the application of a traditional automated theorem prover to resolve the \textsl{Robbins problem} about Boolean algebras \citep{McCune1997}.  The second involves the application of a $\textsc{Sat}$-solver to resolve the \textsl{Boolean Pythagorean Triple conjecture} in arithmetical combinatorics \citep{Heule2016}.  And the third involves the use of a large language model (in concert with other techniques) to improve a lower bound related to the study of \textsl{cap sets} \citep{Romera-Paredes2024}, also a topic from combinatorics.

Popular reporting on these results might lead one to expect that we have already reached the point where automated techniques play a significant role in settling open questions.\footnote{The fourth quote in the epigraph is a classic example.  But see also, e.g., \citep{Roberts2023}, \citep{Sample2023},\citep{Castelvecchi2024}.}  But we will also suggest that they have the following characteristics:
\begin{enumerate}[1)]
\item The techniques by which they have been demonstrated are implementations of \textsl{brute force search} -- i.e. refinements of a generate-and-test procedure.  
\item The results which have been obtained are of \textsl{low logical complexity} -- i.e. they are $\Sigma^0_1$-statements when formalized in the language of first-order arithmetic.
\end{enumerate}   

In \S 5 we will suggest that these observations highlight several under-explored questions about the role of computation in mathematics -- e.g.   What exactly is meant by ``brute force''?  How is a demonstration conducted in this manner different than a traditional mathematical proof?  To what extent does logical complexity track discovermental complexity?  But these considerations notwithstanding, we will conclude by suggesting that characteristics 1) and 2) illustrate why current and foreseeable results obtained via artificial intelligence-inspired techniques do not tell against the inherent difficulty of proof discovery.

\section{The proof discovery problem}

Here is an initial formulation of our main thesis:
\begin{itemize}
\item[(N$_0$)]
Advances in computing theory or technology -- inclusive of those currently understood to exemplify artificial intelligence -- will not radically alter the following aspect of our current understanding of mathematics and its practice:
\textsl{Settling the status of open questions is an inherently difficult problem}.
\end{itemize}
In this section we will clarify this claim relative to what we refer to as the \textsl{proof discovery problem} and explain why we take this to be a central aspect of mathematical practice.  This  requires also clarifying what we mean by ``open questions'' and ``inherent difficulty''.   

By an open question we mean a mathematical statement for which we currently lack a proof or a refutation which we take to establish that the statement in question is true or false.    This notion is readily exemplified by what are often called ``substantial'' or ``non-trivial'' open questions of the sort exemplified by those appearing on lists such as the Millennium Problems \citep{Carlson2006} or those of \citet{Hilbert1900}, \citet{Landau1912}, \citet{Smale1998}, or \citet{Nash2016}.  Such collections include statements with elementary number-theoretic formulations (the existence of infinitely many perfect numbers, the Twin Prime Conjecture, the Goldbach Conjecture) as well as those which are known to possess equivalent formulations in the language of arithmetic (e.g. the Riemann Hypothesis, $\mathbf{P} \neq \mathbf{NP}$, the Collatz or $3n+1$ conjecture).   But they also include statements from several other branches of pure and applied mathematics -- e.g. the Naiver-Stokes existence and smoothness problem about fluid flow, the Birch and Swinnerton-Dyer conjecture about elliptical curves, the Hodge conjecture in algebraic and complex geometry, and the Novikov conjecture about the topology of higher-dimension manifolds.      

Many interesting questions can be raised directly about the history and significance of these specific statements.  But of course they are also examples of a great many other currently open questions which will be known to specialists at a given time.  As such, one feature which we take to be part of the common notion of an open question is that such statements form an indefinitely extendable class.  

This is so in the obvious sense that mathematicians continue to formulate novel conjectures, some of which are quickly resolved, others not.   But it is also true in the perhaps trivial seeming sense that since the language in which we carry out mathematics enables us to express infinitely many distinct propositions, there will always be infinitely many statements for which we currently lack a proof or refutation.   From here it is only a short step to the observation that proof discovery can be regarded as a decision problem with infinitely many instances.  In order to formulate this observation precisely we will now take several steps towards formalization which might initially seem counter to our ultimate goal of arguing that resolving open questions is difficult in a \textsl{formalism-independent} sense.  

The first of these is to assume that the sorts of proofs in which mathematicians actually traffic can, at least in principle, be formalized relative to an accepted set of axioms $\Gamma$ and notion of provability $\vdash_S$ by which theorems are derived from axioms.   The limitative theorems on which our basic argument is based are \textsl{highly invariant} to exactly how these simplifying assumptions are made.  But for purposes of concreteness, it will still be useful to fix what might be regarded as consensus choices for $\Gamma$ and $\vdash_S$.  To avoid becoming embroiled in details, we will thus assume that $\Gamma$ is a computably axiomatized theory over a countable language $\mathcal{L}_{\Gamma}$ sufficient for formalizing contemporary core mathematics -- e.g. $\mathsf{ZFC}$ or an extension with a few large cardinal axioms or a similar system such as the Tarski-Grothendieck set theory will be more than sufficient.   We will similarly take $\vdash_S$ to denote the derivability relation in a conventional proof system for (classical) first-order logic.\footnote{We are fully aware that the assumption that ``all of mathematics'' within a fixed axiomatic system will strike some readers as contentious for a variety of reasons. But it will become clear below that the adoption of anti-formalistic standpoints will tend to buttress rather than detract from the argument given below.}

What we will refer to as the \textsl{proof discovery problem} is thus simply that of determining if a given $\mathcal{L}_{\Gamma}$ sentence $\varphi$ is derivable from $\Gamma$ in the sense of $\vdash_S$. This is in turn equivalent to deciding membership in the set 
\begin{example} $\textsc{Prov}^S_{\Gamma} = \{\varphi \in \mathcal{L}_{\Gamma} : \Gamma \vdash_S \varphi\}$ 
\end{example}
$\textsc{Prov}^S_{\Gamma}$ can thus be compared to other decision problems like primality checking.  This is similarly represented as a set $\textsc{Primes} = \{n \in \mathbb{N} : n \textsl{ is prime}\}$ which gives rise to infinitely many individual yes/no questions of the form $n \in \textsc{Primes}$.

Deciding membership in $\textsc{Primes}$ is a task we typically perform using a primality algorithm -- e.g. trial division or the sieve of Eratosthenes.  By applying such an algorithm, we are able to determine both positive and negative propositions about primality -- e.g. $666889 \in \textsc{Primes}$ but $667071 \not\in \textsc{Primes}$.   Most readers will of course be aware that for many choices of $\Gamma$ and $\vdash_S$,  there can be no algorithm for uniformly deciding if $\varphi \in \textsc{Prov}^S_{\Gamma}$ or $\varphi \not\in \textsc{Prov}^S_{\Gamma}$ (and this will indeed figure in our argument below).  But the point which we wish to stress initially is that the task of deciding membership in this set can be understood as at least \textsl{approximating} the task of resolving open questions in mathematical practice.   

Having now formulated proof discovery as a decision problem, we turn to the task of motivating the first two premises of the argument we will give in \S 3:

\begin{itemize}
\item[(P1)]
\begin{enumerate}[i)]
\item \textsl{Proof discovery} is a central goal of mathematical practice.
\item The \textsl{difficulty} of proof discovery is accurately measured by the classification of the decision problem $\textsc{Prov}^S_{\Gamma}$ with respect to the hierarchies of computability theory and theoretical computer science.   
\end{enumerate}
\end{itemize}

We will take (P1.i) to largely speak for itself.\footnote{The centrality which mathematicians themselves assign to settling open questions was perhaps most famously sloganized by \citet{Hilbert1930a} as ``\textsl{Wir m\"ussen wissen, wir werden wissen}''.   But see, e.g., \citep{Hardy1940} or \citep{Smorynski2020} for similar assessments.} But in highlighting the centrality of deciding open questions to mathematics, we by no means wish to minimize the importance of other aspects of practice such as \textsl{calculation}, \textsl{proof verification}, and \textsl{conjecture discovery}. Computing technology has and continues to be successfully applied to these tasks. The best-known prior examples of calculation and verification have for the most part relied on ``classical'' techniques such as the application of (provably correct) numerical algorithms or formal verification systems based on traditional deductive calculi.\footnote{Examples of the first sort include the hand calculation which underlies Lucas's 1876 demonstration that the 39-digit Mersenne number $M_{127}$ or the the use of computer for case checking in the proof of the Four Color Theorem.  (See, e.g., \citep{Williams1998}, \citealp{Detlefsen1980}, \citealp{Appel1977a}, and \citealp{Tymoczko1979}).  A well-known example of the second sort is the use of a formal verification system or (or \textsl{proof assistant}) in checking Hale's proof of the Kepler Conjecture.  (See, e.g., \citealp{Hales2017} or \citealp{Avigad2021}.)} But there have also been recent successful applications of conjecture discovery using inductive techniques such as machine learning which are more closely associated with artificial intelligence.\footnote{Classical sources for the use of inductive methods in conjecture discovery and the automation thereof of include \citep{Polya1954} and \citep{Borwein2008}.  Applications of machine learning techniques within this frame include \citep{Gauthier2016}, \citep{Davies2021}, and \citep{Gauthier2023a}.}

We mention these tasks to highlight that our thesis is not meant to exclude such technologies coming to play a substantial role in mathematical research.  Indeed it seems all but assured that they will continue to increase in significance.  Rather our argument is intended to illustrate \textsl{that proof discovery is of a fundamentally different nature}.  While this will become clearer as we go on, the basic contrast is  that proof discovery prototypically involves a form of \textsl{unbounded search} which is not required for algorithmic calculation, proof verification, or conjecture discovery.

Our basic contention in regard to (P1.ii) -- which we will clarify further in \S 3 -- is that the classification of the decision problem $\textsc{Prov}^S_{\Gamma}$  with respect to structures like the arithmetical and polynomial hierarchies provides at least some gauge of the difficulty which we confront in resolving open questions in practice -- or as it is sometimes called their \textsl{discovermental complexity}.\footnote{This term was coined by Detlefsen who characterized the discovermental (or \textsl{inventional}) complexity of a mathematical proposition as that which is ``encountered in coming up with a proof in the first place''  \citeyearpar[p. 376]{Detlefsen1990}.  See also  \citep{Detlefsen1996}, \citep{Dean2019c}, and \citep{Arana2023}.}  This may also speak for itself from the perspectives of logic and theoretical computer science.   But to avoid misunderstandings we will first canvas some potential objections.

\subsection{Formalism}  An obvious concern in identifying the task of resolving open  questions with deciding membership in $\textsc{Prov}^S_{\Gamma}$ is that we have adopted a formalistic understanding of mathematics which replaces genuine concern for the truth or falsity of a statement with the question of its derivability from certain axioms.  This problem is illustrated by formally independent statements which will exist for even maximalist choices of $\Gamma$ and $\vdash_S$.   In particular, these represent statements for which both $\varphi \not\in \textsc{Prov}^S_{\Gamma}$ and also $\neg \varphi \not\in \textsc{Prov}^S_{\Gamma}$.  

But of course for at least some such undecidable statements we think that there is an independent fact of the matter about their truth value -- e.g. the consensus view is that $\mathrm{Con}(\mathsf{ZFC})$ is \textsl{true} despite its undecidability in $\mathsf{ZFC}$.\footnote{Such metamathematical statements abou non-derivability often depend on the consistency of the relevant theory $\Gamma$.   But we will assume that this is not contentious here.}  On this basis one might reasonably object to the proposal that deciding membership in $\textsc{Prov}^S_{\Gamma}$ is a reasonable surrogate for resolving open mathematical questions as it is based on a conflation of \textsl{falsity} with \textsl{formal non-derivability}. 

Note, however, that the traditional contention that ``truth outstrips provability''  concerns theories $\Gamma$ which have \textsl{intended models} $\mathcal{M} \models \Gamma$.  But on the other hand, many of the specific domains in which the techniques described in \S 4 have been successfully employed concern so-called ``algebraic'' theories which lack such models.    On the other hand, in the case where $\Gamma$ is a arithmetical or set theoretic theory,  it can often be demonstrated that the class of statements true in $\mathcal{M}$ --  $\mathrm{Th}(\mathcal{M})$ -- is of (much) higher complexity that $\textsc{Prov}^S_{\Gamma}$.   Thus although $\textsc{Prov}^S_{\Gamma} \subset \mathrm{Th}(\mathcal{M})$ (in light of the \textsl{soundness} of $\Gamma$ and $\vdash$ with respect to truth in $\mathcal{M}$), deciding the \textsl{truth} of arbitrary $\mathcal{L}_{\Gamma}$ statements will typically be \textsl{harder} than deciding their provability.\footnote{A canonical illustration is provided by first-order Peano arithmetic $\mathsf{PA}$ whose standard model is $\mathcal{N}$ with natural numbers $\mathbb{N}$ as domain and the standard interpretation of addition and multiplication.  In this case $\textsc{Prov}^S_{\mathsf{PA}}$ is a $\Sigma^0_1$-definable set (for $S$ a standard proof system for first-order logic) whereas $\mathrm{Th}(\mathcal{N})$ is a $\Delta^1_1$-complete set whose complexity is ``properly above'' not only $\textsc{Prov}^S_{\mathsf{PA}}$ but also the entire arithmetical hierarchy.}     

It follows that theorists who wish to object to (P1.ii) in light of their rejection of formalism will also be preinclined to view the problem of deciding open questions as \textsl{more difficult} than this premise prescribes.   But this  sets the bar (considerably) higher for what a computational system would have to achieve in order to surpass what such theorists believe human mathematicians are already able to accomplish -- i.e. deciding mathematical truth as opposed to (mere) provability.   Thus while such theorists may find little interest in the argument rehearsed below, this is so because they are already preinclined to accept its conclusion for other reasons.\footnote{E.g. proponents of anti-mechanist proposals such as that of \citep{Lucas1961}.}

\subsection{Formal versus informal provability}  A related concern is that in equating the difficulty of proof discovery with that of deciding membership in $\textsc{Prov}^S_{\Gamma}$, we distort the kind of challenge which we face in resolving open questions in practice by focusing on \textsl{formal} rather than \textsl{informal} provability.  Such an identification is required by our argument below for two reasons.   First, as we will see in \S 3, the argument is mediated by a premise citing mathematical results about the complexity of $\textsc{Prov}^S_{\Gamma}$ which would not be available if we were to base our characterization of proof discovery on informal provability.  Second, current automated techniques are designed to search either directly for formal proofs or other sorts of finitary objects which witness the truth of a mathematical statement in a manner which can be converted into a formal proof.  

These considerations notwithstanding, the sort of objection to (P1.ii) we have in mind might be leveled by theorists dissatisfied with the view that an adequate demonstration of a mathematical proposition must take the form of a formal proof from declared axioms rather than the sort of informal demonstration of the sort in which we trade in mathematical practice.   On its own, this does not precisely characterize the relevant distinction between formal and informal proofs.  But some dimensions along which the latter are taken to differ from the former include a lack of an explicit choice of formal language or axioms, reliance on semantically mediated inferences rather than explicitly stated rules of inference, and compression of steps in favor of reliance on ``high-level outlines'' or diagrams.\footnote{For instance \citet{Avigad2021} states that ``an informal proof is a form of data compression'' (p. 7386).  A simple example discussed by \citet[p. 4]{Wiedijk2006} involves the familiar informal proof of the irrationality of $\sqrt{2}$ wherein we assume for reductio that there exist $a$ and $b$ such that $a^{2} = 2b^{2}$ with greatest common divisor $1$. This formulation implicitly assumes that the fraction $\frac{a}{b}$ is given in lowest terms.   When this proof is formalized in a manner such that it can be verified by an automated verification system, it must be made explicit that this may be assumed without loss of generality.  The process of making such assumptions explicit -- either by treating them as axioms or supplying other formalized inductive proofs -- may thus indeed be described as ``decompressing''  information which is ``hidden'' (or implicit) in the original informal proof.}

The example we will consider in \S 4.2 illustrates an instance in which a formal proof is the \textsl{only} sort of demonstration which may be available to resolve an open question.  But the more common scenario is that a proposed resolution of an open problem will not be announced unless it is judged by its authors to contain information sufficient for the  mathematical community to evaluate its correctness. Such claimed demonstrations are then scrutinized as part of the normal refereeing procedure. In a small minority of cases this results in the judgement that formal verification is required to assess either the truth of the statement in question or the validity of the proposed informal proof. 

Experience has borne out that this latter step can  be arduous, potentially requiring much effort by human mathematicians skilled in the use of verification systems and potentially also substantial computer time to fill in intervening steps.  This in turn illustrates that there may indeed be a distinct sense in which finding an initial informal proof of an open question $\varphi$ is typically \textsl{easier} than constructing a formal derivation from $\Gamma$ which witnesses $\varphi \in \textsc{Prov}^S_{\Gamma}$.   There is thus a concern that the complexity of deciding $\varphi \in \textsc{Prov}^S_{\Gamma}$ will \textsl{overestimate} the difficulty which we face in settling the status of $\varphi$ in the sense we care about in practice.

We will not contest this observation as such.  But we propose that a principled distinction may still be made between the sort of difficulty involved in originally finding an informal proof and that of transforming such a proof into a formal derivation.   For experience also bears out that the initial difficulty of finding a ``high level'' argument which succeeds in convincing the mathematical community -- say on the basis of published text $\mathcal{D}_1$ -- still typically overshadows that of converting such a proof into a formal derivation $\mathcal{D}_2$  from a given set of axioms $\Gamma$.  In particular, while finding $\mathcal{D}_1$ at least sometimes involves what is subsequently described as the ``novel'', ``creative'', or ``ingenious'' discovery of new ideas and techniques, constructing $\mathcal{D}_2$ from $\mathcal{D}_1$ typically involves only ``meticulous'' or ``mechanical'' attention to ``routine'' details such as finding appropriate definitions and lemmas to fill in gaps, proving the correctness of the implementations of algorithms, etc.\footnote{\citet{Avigad2021} provides a detailed account of how and why formalization helps to underpin the reliability of informal proof.  On the specific process of finding formal proofs from informal ones he writes ``When someone $\ldots$ embarks on a formalization project, the assumption that the theorem can be formalized is never in question $\ldots$ The question is, rather, how best to go about it and how long it will take. To be sure, the formalization process sometimes uncovers minor errors and omissions in the informal presentation that have to be remedied, but these are conventionally viewed as just that $\ldots$ rather than indications that the informal source is correct but unformalizable.'' (p. 7339)}  

The formal classification of the complexity of $\textsc{Prov}^S_{\Gamma}$ does indeed compound these two forms of difficulty.  But it still seems reasonable to assume that the difficulty of finding $\mathcal{D}_1$ in the first place will typically dominate over that of finding $\mathcal{D}_2$ from $\mathcal{D}_1$.  It is for this reason that we believe that the composite measure still leads to an informative form of the argument given below.  

\subsection{Provability versus proofs} Another potential objection to (P1.ii) draws attention to the fact that while the complexity-theoretic classification of $\textsc{Prov}^S_{\Gamma}$ concerns the infinite class of formal theorems of $\Gamma$, in practice we are concerned with specific instances of the question $\Gamma \vdash_S \varphi$ in cases where $\varphi$ formalizes a recognized open question.  This leaves open the possibility that while deciding membership in $\Gamma \vdash_S \varphi$ is hard \textsl{in the general case},   deciding the \textsl{particular instances} in which we are interested in practice -- e.g. where $\varphi$ is the Twin Prime Conjecture or the Riemann Hypothesis -- may be either (relatively) easier or (relatively) harder.\footnote{A  related objection can be formulated in virtue of the existence of statements about membership in  $\textsc{Prov}^S_{\Gamma}$ which can themselves be formulated as mathematical states which admit to ``meta-proofs'' -- e.g.  statements of the form $\varphi \in \textsc{Prov}^S_{\Gamma}$ will typically be expressible by $\Sigma^0_1$-statements of the form $\exists x \mathrm{Proof}^S_{\Gamma}(x, \ulcorner \varphi \urcorner)$ for an appropriate formalized proof predicate which they themselves attempt to prove or refute.   In this case it is sometimes possible to construct statements $\varphi$ such that $\exists x \mathrm{Proof}^S_{\Gamma}(x, \ulcorner \varphi \urcorner)$ has a ``short'' proof from $\Gamma$, any proof of $\varphi$ must be ``astronomically long'' in comparison (see, e.g., Theorem 1.3 of \citealp{Parikh1971}).   But not only are the statements which exemplify this property of a 
contrived metamathematical character, the constructions which yield them also testify to the  fact that they are \textsl{true} (albeit for reasons which can only be formalized in a theory stronger than $\Gamma$ itself).}   

This concern can be illustrated by a comparison to decision problems which are known to admit ``special cases''.  For instance, while deciding a general instance of the primality problem $n \in \textsc{Prime}$ remains (relatively) difficult, deciding whether a \textsl{Mersenne number} -- i.e. one of the form $n = 2^m-1$ -- is prime has long been known to admit to an easier method known as the \textsl{Lucas-Lehmer test}.  In this case the set $\textsc{Mersenne} = \{n \in \mathbb{N} : \exists m(n = 2^m -1) \ \& \ n \text{ is prime}\}$ might be described as an \textsl{easier subproblem} of $\textsc{Primes}$.\footnote{See, e.g., \citep{Crandall2005}.    Of course it is now known that $\textsc{Primes}$ itself admits a polynomial time decision procedure.  Better examples of the phenomenon in question are thus provided by decision problems  (e.g. $\textsc{Sat}$) which are provably hard for a given class (e.g. $\mathbf{NP}$) but also contain precisely delimited subclasses which can be proven to be easier (e.g. $\textsc{2-Sat}$).}    

The analogous concern about mathematical provability can thus be put as follows: Might not the sort of open questions which mathematicians choose to focus on in practice -- i.e. what could be called \textsl{naturally occurring open questions} -- turn out to be either easier or harder to resolve than deciding arbitrary instances of $\textsc{Prov}^S_{\Gamma}$ ``on average''?   If so, can we formally characterize statements in this class or account for what makes them easier or harder to resolve than ``typical'' statements?

Unlike the prior objections, we take this issue to be of considerable import. For it draws attention to the possibility that our practices of electing to study certain mathematical statements but not others tracks a measurable characteristic of their discovermental complexity.   In fact it is just this sort of issue about the contours of mathematical knowledge which recent work in automated theorem proving and artificial intelligence may help to bring into focus.  We will return this issue in \S 5.

\section{The basic argument}

The story of artificial intelligence is often presented as beginning in the late 1950s with the work of figures such as Newell \& Simon, Wang, Davis, and Putnam which we will discuss in the next section.  Already at this time, the prospects for applying computing technology to proof discovery in mathematics were widely discussed.  It was in the context of these developments that the optimistic predictions illustrated by the second and third quotes in the epigraph began to be announced.

Progress was slow throughout the 1970s, during which the phenomenon of $\mathbf{NP}$-completeness was discovered.  This is now regarded as a characterization of an \textsl{intractable} (or \textsl{inherently difficult})  \textsl{problem} -- a concept which had come into greater focus during the 1960s along with the  availability of digital computers and the desire to apply them to problems of practical import.  This class includes many classic examples about planning and search which originated in artificial intelligence.  It was within this context that the phrase ``combinatorial explosion'' was coined by \citet{Lighthill1973} in the course of presenting a form of the argument -- which was then formulated more incisively by \citet{Rabin1974} -- we are about to consider.\footnote{\citet{Lighthill1973} originally defined a ``combinatorial explosion'' as the property of a ``large knowledge base which results from the explosive growth of any combinatorial expression, representing numbers of possible ways of grouping elements of the knowledge base according to particular rules, as the base's size increases''.   He went on to note that the original optimistic hopes about automated theorem proving were ``disappointed through the power of the combinatorial explosion in rapidly cancelling out any advantages from increase in computer power''.}  

These observations contributed to the first of the so-called ``AI winters''.\footnote{See, e.g., \citep[\S 1.3]{Russell2021}.}   We are evidently living in the midst of a subsequent summer.   Nonetheless, our contention is that the following updated form of the Lighthill-Rabin argument remains relevant for the prospects of artificial intelligence in mathematics:

\begin{itemize}
\item[(P1)] i) \textsl{Proof discovery} is a central goal of mathematical practice.

\vspace{.5ex}

\noindent ii) The \textsl{difficulty} of proof discovery is accurately measured by the classification of the decision problem $\textsc{Prov}^S_{\Gamma}$ with respect to the  hierarchies of computability theory and computational complexity theory.      

\vspace{1ex}

\item[(P2)] For a wide range of relevant choices of $\Gamma$ and $\vdash_{S}$, the problem $\textsc{Prov}^S_{\Gamma}$ is of \textsl{high computational complexity} in the sense of (P1.ii).

\vspace{1ex}

\item[(P3)] If a problem is of high computational complexity in the sense of (P2), then it is \textsl{inherently difficult} (or \textsl{intractible}) to decide arbitrary instances of it using computing hardware which we can construct and apply in practice.     
\end{itemize}

By combining P1, P2, and P3 we reach a restatement of our original thesis:
\begin{itemize}
\item[(N$_1$)] The application of artificial intelligence will not lead to a substantial reduction in the difficulty of one of the central goals of mathematical practice.
\end{itemize}
We have already attempted to clarify the relevant notion of ``proof discovery'' in the course of arguing for (P1).  Once the relevant characterizations of ``high computational complexity'', ``intractible'', and ``computing hardware which we can construct and apply in practice'' are in place, we hope that premises (P2) and (P3) will speak for themselves for readers logic and computer science.   But as these notions are again not widely appreciated within philosophy, we will spend the rest of this section clarifying them in the course of arguing explicitly for these premises.

\subsection{Considerations from computability theory}

An initial observation underlying (P2) is that for the most evident choices of $\Gamma$ and $\vdash_S$, the decision problem $\textsc{Prov}^S_{\Gamma}$ is formally undecidable.   This is so because of familiar considerations surrounding G\"odel's First Incompleteness Theorem and Church and Turing's negative solution to the original \textsl{Enstscheindungsproblem} for first-order logic.  But it will still be useful to collect together several related extensions as follows:\footnote{These results \ref{g1} go back essentially to \citep{Godel1931a}, \citep{Rosser1936}, and \citep{Tarski1953}.  Textbook formulations can be found in \citep{Shoenfield1967} or \citep{Cooper2004}.}
\begin{theorem} 
\label{g1}
Suppose that \textnormal{a)} $\vdash_S$ is a computably enumerable derivability relation which extends that of classical first-order logic $\mathrm{[FOL]}$ and \textnormal{b)} $\Gamma$ is a consistent, computably axiomatizable theory that interprets $\mathsf{Q}$ $($i.e. Robinson arithmetic$)$ and which is additionally $\Sigma^0_1$-sound relative to such an interpretation.\footnote{Recall that an $\mathcal{L}_a$-theory $\mathsf{T}$ is $\Sigma^0_1$-sound just in case  the following condition holds:  for all $\mathcal{L}_a$ statements  of the form $\exists x \varphi(x)$ where $\varphi(x)$ only contains bounded quantifiers, if $\mathsf{T} \vdash \exists x \varphi(x)$, then $\exists x \varphi(x)$ is true in the standard model of arithmetic.  This is equivalent to the condition known as 1-\textsl{consistency} and is satisfied by all the theories mentioned in note \ref{theoriesnote} under the relevant interpretations of $\mathcal{L}_a$ in their languages.   It is necessary here only for part iii) so as to ensure that $\Gamma$ does not prove any false $\Sigma^0_1$-statements which misrepresent the existence of halting computations so as to contrive that the axioms of $\mathsf{T}$ induce a set of theorems which fails to be $\Sigma^0_1$-complete.}   Then:
\begin{itemize}
\item[i)] $\Gamma$ is incomplete relative to the definition of derivability $\vdash_S$.  In particular, there will exist $\mathcal{L}_{\Gamma}$-sentences $\varphi$ such that $\Gamma \not\vdash \varphi$ and $\Gamma \not\vdash \neg \varphi$.   Furthermore, it is possible to find such $\varphi$ which are provably equivalent $\Pi^0_1$-statements in the language $\mathcal{L}_a$ of first-order arithmetic.\footnote{I.e. of the form $\forall x \psi(x)$ for $\psi(x)$ containing only bounded arithmetical quantifiers, either natively or by interpretation into $\mathcal{L}_a$.}
\item[ii)] $\textsc{Prov}^S_{\Gamma}$ is not $\Delta^0_1$-definable and thus formally undecidable as a decision problem.
\item[iii)] $\textsc{Prov}^S_{\Gamma}$ is $\Sigma^0_1$-complete relative to many-one reductions.
\end{itemize}
\end{theorem}

For the sort of consensus choices for $\Gamma$ and $\vdash_S$ for formalizing core mathematics envisioned in \S 2,  conditions a) - b) will typically satisfied in a paradigmatic manner.\footnote{E.g (and \textsl{mutatis muntandis}) if $\Gamma = \mathsf{PA}, \ldots, \mathsf{Z}_2$ (second-order arithmetic), $,\ldots, \mathsf{Z}$ (Zermelo set theory), $,\ldots,  \mathsf{ZF(C)}, \ldots$ and the vast majority of their familiar high-order or class-based extensions, variants based on intuitionistic logic, type-theoretic or categorical reformulations, etc. \label{theoriesnote}}
%\footnote{For instance ``classical foundational theories'' based on their arithmetical or set-theoretic signatures ars more than sufficient for interpreting $\mathsf{Q}$ and are both computably axiomatizable and are $\Sigma^0_1$-sound relative to the relevant interpretations.    This includes first and second-order Peano arithmetic ($\mathsf{PA}$ and $\mathsf{Z}_2$), first-order Zermelo or Zermelo-Fraenkel set theory (i.e. $\mathsf{Z}$ and $\mathsf{ZF}$), classed-based alternatives (e.g. $\mathsf{GB}$ and $\mathsf{MK}$), and the extension of such theoreis with the Axiom of Choice, computable classes of large cardinal axioms, etc.  Similar facts hold for familiar ``constructive foundational theories'' -- e.g. $\mathrm{HA}$, $\mathrm{HA}^{\omega}$, $\mathsf{ITT}$, $\mathsf{HoTT}$, etc.  In the case of theories based on intuitionistic logic it becomes more complicated to relate formula complexity to computational complexity (as we do below) in light of the failure of the Prenex Normal Form Theorem.  Convinced intuitionists may thus have a  different perspective on how computability theoretic-results bear on genuine computational difficulty.  But these concerns will not typically complicate their appreciation of complexity-theoretic results since, e.g., the non-collapse of the Polynomial Hierarchy -- if true -- is compatible with even strong constructive principle like the intuitonistic form of Church's Thesis (i.e. the  formalization of \textsl{all functions are recursive}).   \label{theoriesnote}}  
It thus follows that for a wide range of relevant mathematical theories and notions of derivability, the proof discovery problem is \textsl{as hard as its metamathematical definition allows} in the sense measured by the arithmetical hierarchy.  This can be further glossed as follows:
\begin{enumerate}[(P2.1)]
\item Conditions a) and b) are sufficient to ensure that $\textsc{Prov}^S_{\Gamma}$ is in the class of problems which can be defined as a  $\Sigma^0_1$-formula of the language $\mathcal{L}_a$ of first-order arithmetic.   In particular we will have $$\textsc{Prov}^S_{\Gamma} = \{\ulcorner \varphi \urcorner \in \mathbb{N} : \exists x \mathrm{Proof}^S_{\Gamma}(x, \ulcorner \varphi \urcorner)\}$$ for a suitable arithmetized proof predicate involving only bounded quantifiers and G\"odel numbering $\ulcorner \cdot \urcorner$ of $\mathcal{L}_{\Gamma}$.   
\item On the other hand, part ii) of Theorem \ref{g1} reports that this classification cannot be improved to $\Delta^0_1$-definability -- i.e. the prior classification is \textsl{optimal} with respect to the arithmetical hierarchy.
\item If we consider the set $A = \{\ulcorner \varphi \urcorner : \varphi \in \textsc{Prov}^S_{\Gamma}\} \subseteq \mathbb{N}$, then every other $\Sigma^0_1$-definable set $B \subseteq \mathbb{N}$ is definable from $A$ as a parameter in the form $B = \{n : \psi(n,A)\}$ where $\psi(x,X)$ is a $\Delta^0_1$-formula.   
\end{enumerate}

The foregoing facts pertain to the classification of the \textsl{logical complexity} of $\textsc{Prov}^S_{\Gamma}$ -- i.e. how complicated a formula is required to define it as a set relative to the standard model of arithmetic.  But as many readers will be aware, (P1.1-3) respectively correspond  to properties concerning the \textsl{computational complexity} of $\textsc{Prov}^S_{\Gamma}$  -- i.e. how hard is it decide by an algorithm.  In particular, if $\Gamma$ satisfies the hypotheses of Theorem 1, then $\textsc{Prov}^L_{\Gamma}$ will be a \textsl{many-one complete set} -- i.e.
\begin{enumerate}[(P2.1$'$)]
\item  $\textsc{Prov}^S_{\Gamma}$ is a \textsl{computably enumerable} (or \textsl{semi-decidable}) set -- i.e. there exists a Turing machine $M(x)$ such that if $\Gamma \vdash_S \varphi$, then $M(\ulcorner \varphi \urcorner)$ will  eventually halt outputting ``yes'' and fail to halt otherwise.  
\item $\textsc{Prov}^S_{\Gamma}$ is not a \textsl{computable} (or \textsl{decidable}) set -- i.e. there is no Turing machine $M(x)$ such that for all $\varphi \in \mathcal{L}_{\Gamma}$, $M(\ulcorner \varphi \urcorner)$ halts outputing ``yes'' if $\Gamma \vdash_S \varphi$ and if  $\Gamma \not\vdash_S \varphi$, then  $M(\ulcorner \varphi \urcorner)$ halts outputting ``no''.
\item  The existence of a decision algorithm for $\textsc{Prov}^S_{\Gamma}$ would imply a decision algorithm for \textsl{every} $\Sigma^0_1$-set $B$ in the following sense: there exists a Turing computable function $f(x)$ such that for all $n \in \mathbb{N}$, 
\begin{center}
$n \in B$ iff $f(n)$ is the G\"odel number of a $\mathcal{L}_{\Gamma}$ formula $\varphi$ such that $\Gamma \vdash_S \varphi$.
\end{center}
\end{enumerate}

The existence of such correspondences between descriptive and computational complexity stands behind the use of the expression \textsl{degree of difficulty} to refer to levels of arithmetical hierarchy and related structures in computability theory.\footnote{This expression originates with \citep{Post1944}.  See, e.g., \citep{Slaman1998} for more on using formula complexity and the related definition of the Turing degrees to measure mathematical difficulty.}  Per premise (P2) of the basic argument, the foregoing results are illustrative of the ``high computational complexity'' of $\textsc{Prov}^S_{\Gamma}$  in the case that $\Gamma$ and $\vdash_S$ satisfy the property a) and b) of Theorem \ref{g1}.  Of course the term ``high'' should be understood here as relative to the use of the arithmetical hierarchy as a scale for measuring some salient notion of difficulty related to mathematical practice.\footnote{In particular, a $\Sigma^0_1$-complete set like $\textsc{Prov}^S_{\Gamma}$ is only of ``high complexity'' in the sense of being ``above'' (in the sense of reducibility) sets which are either computable (i.e. $\Delta^0_1$) or computably enumerable but not $\Sigma^0_1$-complete.   For instance $\textsc{Prov}^S_{\Gamma}$ is ``below'' -- and thus in the relevant sense \textsl{less difficult than} -- a $\Pi^0_2$-complete set such as $\textsc{Tot}$ -- i.e. set of indices to Turing machines which halt on all inputs.  Two relevant aspects of the use of ``high'' in the formulation of (P2) are thus as follows: i) this complexity is \textsl{maximal} given that $\textsc{Prov}^S_{\Gamma}$ admits a $\Sigma^0_1$-definition; ii) a many-one complete set can be regarded as \textsl{astronomically more difficult} to decide in the general case than the decidable sets classified by the complexity-theoretic hierarchies considered below.} But in the case of comparing $\textsc{Prov}^S_{\Gamma}$ with a decidable problem like $\textsc{Primes}$, the following sub-argument can also now be adduced in favor of (P3):\footnote{Note that this argument avoids recourse to Church's Thesis.  In the relevant case, this would entail that since $\textsc{Prov}^S_{\Gamma}$ is not Turing-decidable, it is also not decidable by an informally specified algorithm.  But the concern here is not with whether humans might be able to exceed computers by carrying out algorithms which are intuitively effective but not Turing computable.  Rather it is whether the sort of concrete computing systems we can construct and apply in practice might be of assistance in deciding instances of $\textsc{Prov}^S_{\Gamma}$ which confound human mathematicians.}
\begin{enumerate}[(P3.1)]
\item Since by P2.2$'$) $\textsc{Prov}^S_{\Gamma}$ is not decidable by a Turing machine, it follows that it is also not decidable by an abstract machine from any \textsl{model of computation} which computes the same class of functions as the basic Turing machine model $\mathfrak{T}$ in terms of which computability theory is standardly developed.   
\item Not only does this class include the familiar models of the 1930s --  general recursive functions, the untyped lambda-calculus, etc. -- it also includes a broad class of what might be called \textsl{generalized finitary} models which satisfy the \textsl{locality} and \textsl{boundedness} conditions originally formulated by \citet{Gandy1980} and refined by \citet{Sieg2009}.   This class includes Gandy machines and cellular automata, but also the abstract characterizations of the architectures of contemporary computers -- e.g. the \textsl{random access machine} model $\mathfrak{R}$ and many of its variants.
\item As such, contemporary computing machinery does not allow us to uniformly decide instances of $\textsc{Prov}^S_{\Gamma}$ in the following sense:  there does not exist a concrete computing device -- i.e. what we would now conventionally call a ``computer'' -- which we can both construct and employ in practice to concretely decide arbitrary instances of $\textsc{Prov}^S_{\Gamma}$.
\item In the case where a theory $\Gamma$ satisfies the hypotheses of Theorem \ref{g1}, another form of confirmation of (P3.3) derives from the $\Sigma^0_1$-completeness of $\textsc{Prov}^S_{\Gamma}$.  For as reported by (P2.3$'$), a decision method for this problem would lead to such a method for a number of other $\Sigma^0_1$-definable problems of both theoretical and practical import -- e.g. the Halting Problem, Hilbert's Tenth Problem, the word problem for semigroups, etc.   These are illustrative of problems which experience suggests we also cannot decide in practice.
\end{enumerate}

\subsection{Considerations from complexity theory} The foregoing points are familiar and were appreciated early on.\footnote{In fact the import of the basic argument appears to have been largely foreseen at the time of the original formulation of the \textsl{Entcheindungsproblems} -- i.e. does there exist an algorithm for determining whether a first-order formula is valid (and hence provable from no premises)?   In their original textbook statement, \citet[p. 77]{Hilbert1928} describe this as ``the main problem of mathematial logic''.  For as they observe, if $\Gamma$ is a finitely axiomatizable theory  the problems can be understood as  equivalent in the sense that $\varphi \in \mathrm{Prov}^{\mathrm{FOL}}_{\Gamma}$ iff $\bigwedge \Gamma \rightarrow \varphi \in  \mathrm{Prov}^{\mathrm{FOL}}_{\emptyset}$.  This includes not only the axiomatic theories of geometry which Hilbert and Ackermann used to movitate the problem, but also $\mathsf{Q}$ itself as well as far stronger theories like G\"odel-Bernays set theory.  From this it can  be taken to follow that if the   \textsl{Entcheindungsproblem} had a positive solution, then the general proof discovery problem in mathematics could be solved algorithmically.  Hilbert and his collaborators were momentarily optimistic about this in light of positive results for restricted formula classes.  But per the first quote in the epigraph, \citet{Hardy1929} (and even earlier \citealp{Neumann1927}) expressed skepticism as to whether a decision algorithm would be forthcoming and also realized the ramifications a positive resolution would have.}  But the potential applications of artificial intelligence to proof discovery require us to proceed more carefully in arguing for (P3).  This is so in light of at least two concerns.   First, at least one of the applications we will consider in \S 4 pertains to choices of $\Gamma$ and $\vdash_S$ which do not directly fall within the scope of Theorem \ref{g1}.   Second, one might also think advances in computing technology which will themselves be brought about by artificial intelligence will make the sub-argument (P3.1-4) obsolete.  

In order to appreciate the first point, recall that \textsl{Presburger arithmetic} $[\mathsf{PrA}]$ can be characterized as the set of sentences in the language $\mathcal{L}^-_a$ consisting of $\mathcal{L}_a$ without the multiplication symbol which are true in in the structure $\mathcal{Z} = \langle \mathbb{Z},<,+,0,1 \rangle$. While this is a semantic definition, \citet{Presburger1930} originally presented an axiomatization of $\mathsf{PrA}$ consisting of the axioms of an ordered abelian group together along with a schema stating the basic property of Euclidean division. Using the notation we have adopted, it can then be shown that $\textsc{Prov}^{\mathrm{FOL}}_{\mathsf{PrA}} = \mathrm{Th}(\mathcal{Z})$.\footnote{See \citep[III.4]{Smorynski1991a} for a modern formulation.}

Since $\mathsf{Q}$ cannot be interpreted in $\mathsf{PrA}$, it follows that Theorem \ref{g1} does not apply.   In fact \citet{Presburger1930} originally showed that $\mathsf{PrA}$ is \textsl{decidable} by showing that it admits quantifier elimination.  This procedure then became a natural target to implement on early digital computers starting in the 1950s.  But it was quickly realized that both Presburger's original algorithm and several refinements could not be carried out concretely even for very simple formulas.\footnote{This was already observed by \citet{Davis1957}. See \citep{Haase2018} for a survey of later developments.}  

This provided the context for Michael Rabin's \citeyear{Rabin1974}  paper ``Theoretical impediments to artificial intelligence''  whose point of departure was the following result:\footnote{As reported in \citep{Fischer1974}, Rabin obtained this result with Fischer soon after the announcement of Theorem \ref{cl} which is traditionally understood to have initiated computational complexity theory in its contemporary form.}

\begin{theorem} \label{fr} There exists a constant $c > 0$ such that for every  algorithm $\alpha$ deciding $\mathsf{PrA}$, there exists an integer $n_0$ such that for every $n > n_0$ there exists a $\mathcal{L}^-_a$-sentence $\varphi$ of length for which $\alpha$ requires more than $2^{2^{cn}}$ computational steps to decide if $\varphi \in \textsc{Prov}^{\mathrm{FOL}}_{\mathsf{PrA}}$.
\end{theorem}

Rabin's original presentation makes clear that the class of algorithms in question can be taken to encompass not only those implementable relative to standard models of computation such as $\mathfrak{T}$ or $\mathfrak{R}$ but also their \textsl{non-deterministic} variants.\footnote{In fact \citet{Berman1980} subsequently showed that $\textsc{Prov}^{\mathrm{FOL}}_{\mathsf{PrA}}$ is \textsl{complete} for a class known as $\mathbf{STA}(*,2^{2^{cn}},n)$ based on an alternating model which extends doubly exponential nondeterministic time (\textbf{2-NEXP}) and thus also $\mathbf{PSPACE}$ .}  Thus while $\mathsf{PrA}$ is decidable in principle, Rabin suggested that both this result and our practical experiences bear out the fact that it is undecidable in practice for concrete formulas of surveyable length.  

Taking this in conjunction with other early results about the computational complexity, he summarizes the situation as follows:

\begin{quote}
\footnotesize{The above results $\ldots$ point to the possibility that this bad behavior of theorem proving programs is inherent. That [is to say] at least as long as we have a general purpose theorem prover either for first-order logic, or for the whole of propositional logic, or for all statements of even a very elementary and decidable fraction of mathematics, the following bad thing will happen. For rather short sentences it will very often, perhaps even almost always, be the case that the shortest proof or shortest decision algorithm will be impractically long to an extent that settling the question whether the sentence is a theorem, or producing a proof for the sentence, will be impossible in practice. \\
\hspace*{1ex} \hfill \citep[p. 617]{Rabin1974}}
\end{quote} 
Relative to the basic argument, Rabin's moral may thus be framed as follows: even if we restrict attention to cases in which $\textsc{Prov}^S_{\Gamma}$ is decidable in the in-principle sense of computability theory, this problem still often turns out to be \textsl{intractable} (or \textsl{inherently difficult}) to decide in the in-practice sense of computational complexity.\footnote{\citet{Fischer1974} also show that the analogue to Theorem \ref{fr} also holds with respect to any axiomatization $\mathsf{T}$ of $\mathrm{Th}(\mathcal{Z})$ in which axiomhood is decidable in polynomial time -- i.e. one such that $\textsc{Prov}^{\mathrm{FOL}}_{\mathsf{T}} = \mathrm{Th}(\mathcal{Z})$ -- then there are true statements $\varphi \in \mathrm{Th}(\mathcal{Z})$ whose shortest proof in $\mathsf{T}$ is also doubly exponential in the size of $\varphi$.  In regard to his titular concern, Rabin goes on to predict that results in this vein inevitably impose limitations on subsequent developments in artificial intelligence in light of the fact that many reasoning tasks will involve logical deduction.  And in fact many practical illustrations of this pattern many now be cited (see e.g. \citealp[\S 11, \S 13, \S 16]{Russell2021}).}  

Of course Theorem \ref{fr} itself bears directly on the difficulty of proof discovery in mathematics at large, we would also need to inquire further into whether the sort of naturally occurring open questions we considered in \S 2 are expressible in a restricted language such as that of $\mathcal{L}^-_a$.\footnote{An evident concern in this case is that without a means of defining multiplication uniformly in the language, there is no way of expressing open number theoretic conjectures such as the Goldbach or Twin Prime conjectures. A potentially better example to make the current point is those the theory of  \textsl{Real Closed Fields} $[\mathsf{RCF}]$ in whose language arbitrary statements of Euclidean geometry can be stated.  \cite{Tarski1959a} famously showed that $\mathsf{RCF}$ is decidable.  But \citet{Fischer1974} showed a version of Theorem 2 also holds for $\textsc{Prov}^{\mathrm{FOL}}_{\mathsf{RCF}}$ with the time bound $2^{cn}$ -- i.e. while decidable in principle, the proof discovery problem for $\mathsf{RCF}$ is is still exponentially hard and thus intractable from the standpoint of contemporary complexity theory.}  In order to understand the significance of this, it will be useful to state and comment briefly on a yet better known result:     

\begin{theorem} \label{cl}
The set $\textsc{Sat}$ of satisfiable formulas of classical propositional logic is complete for the complexity class $\mathbf{NP}$ consisting of problems decidable in time $O(n^k)$ by a non-determistic Turing machine.     
\end{theorem}

This is to say that $\textsc{Sat}$ is decidable by a non-deterministic algorithm with running time polynomially proportional to the size of the input formula (measured by the numbers of symbols) and also that every problem decidable in this manner may be efficiently reduced to $\textsc{Sat}$ -- i.e. in a manner analogous to (P2.3$'$) with $f(x)$ itself computable in polynomial time.  Since a sentence  $\varphi$ is a member of the set $\textsc{Taut}$ of propositional tautologies just in case $\neg \varphi \not\in \textsc{Sat}$, this means that  $\textsc{Taut}$ is a so-called $\mathbf{coNP}$-complete problem.  But since $\textsc{Taut}$ also coincides with the set $\mathrm{Prov}^{\mathrm{PL}}_{\emptyset}$ of theorems of classical propositional logic, this means that proof discovery for this system is itself a $\mathbf{coNP}$-complete problem.

Theorem \ref{cl} -- which is now known as the \textsl{Cook-Levin Theorem}\footnote{See, e.g., \citep{Aaronson2016} or \citep{Dean2019c} for discussion of its history and signficance.} -- had only just been formulated at the time of the lecture on which Rabin's paper is based.  But already at this point, a wide range of other combinatorial problems which were (and still are) believed to be intractable in the general case had already been shown to be $\mathbf{NP}$-complete.  This helped establish $\mathbf{NP}$-hardness as a sufficient condition to regard a problem as intractable.   Thus the fact that $\mathrm{Prov}^{\mathrm{PL}}_{\emptyset}$ is complete for the presumptively \textsl{harder} class $\mathbf{coNP}$ can thus be taken to testify that while it is easier to decide than $\textsc{Prov}^{\mathrm{FOL}}_{\mathsf{PrA}}$, it is still an inherently difficult problem in the sense measured by the hierarchies studied in computational complexity theory.\footnote{The situation is complicated by the fact that we do not currently know whether $\mathbf{NP}$ is distinct from the class $\mathbf{P}$ of problems decided in polynomial time by a \textsl{deterministic} Turing machine (which is conventionally taken to coincide with the problems which we can feasibly decide in practice).  It is, however, provable that $\mathbf{NP}$ is a proper subset of the class $\mathbf{EXP}$ of problems decidable in exponential time.  It is in this sense that $\mathrm{Prov}^{\mathrm{PL}}_{\emptyset}$  is (provably) less difficult than $\textsc{Prov}^{\mathrm{FOL}}_{\mathsf{PrA}}$, which is in turn (provably) less difficult than  $\textsc{Prov}^{\mathrm{FOL}}_{\emptyset}$ or $\textsc{Prov}^{\mathrm{FOL}}_{\mathsf{Q}}$, etc.}

Theorem \ref{cl} is relevant to our current concerns in virtue of ongoing attempts we will discuss in \S 4.2 to employ decision algorithms for propositional satisfiability to resolve open mathematical questions.  On the other hand, if the widely believed conjecture that $\mathbf{NP} \neq \mathbf{coNP}$ is true, then there exist infinitely many theorems of propositional logic whose shortest proofs in any conventional proof system are ``infeasibly long'' relative to their lengths, just as in the case for $\mathsf{PrA}$.\footnote{This includes definition of $\vdash_S$ based on Hilbert systems, natural deduction, and sequent calculi.  See \citep{Cook1979} for the relevant definition of a ``conventional proof system''.}

\subsection{On reasonable models and artificial intelligence}

Before concluding our defense of (P3) -- and thus of the basic argument itself -- one other caveat should be noted.  When we transition from understanding ``difficulty'' from the perspective of computability theory to that of complexity theory, more care also needs to be taken in regard to the models of computation to which the (putatively) limitative theorems apply.  For to read off concrete consequences about what we can and cannot prove in practice from Theorems \ref{fr} and \ref{cl} we need to restrict attention to the narrower class of \textsl{reasonable models of computation} in order to employ a complexity-theoretic analogue of (P3.3).   

The notion of such a model plays a heuristic role in complexity theory similar to that played by an \textsl{effective procedure} in the early history of computability theory.  Just as computability theorists have come to accept the equation of the class of problems decidable by effective procedures with those computable by a Turing machine (i.e. Church's Thesis), complexity theorists have come to accept a similar equation of the class of \textsl{feasibly} (or \textsl{practically}) decidable problems with those in the class $\mathbf{P}$ decidable in polynomial time by a deterministic Turing machine with a suitably efficient representation of inputs and outputs (i.e. the Cobham-Edmond's Thesis).  This leads to the following characterization of a reasonable model in the form of what \citet{Boas1990} originally called  the \textsl{Invariance Thesis}: a model $\mathfrak{M}$ is \textsl{reasonable} just in case for every machine $M \in \mathfrak{M}$ there exists $T \in \mathfrak{T}$ which simulates $M$ with polynomial time overhead and constant space overhead.

The class of reasonable models is again robust and includes the RAM-like models on which the architecture of contemporary digital computers are based.\footnote{See, e.g., \citep{Hennessy2003} on concrete implementational details.}   There are, however, generalized finitary models (in the sense of P3.2) which are not reasonable.   A paradigm example of is the so-called  \textsl{Parallel RAM} model $\mathfrak{P}$ which allows a single processor to spawn a fixed finite number of successor processors sharing a common memory at each step.  This model has been of theoretical interest in complexity theory as it facilities the study of problems which can be efficiently decided by parallelism of \textsl{bounded depth}.\footnote{See, e.g., \citep{Greenlaw1995}.} 

Using parallelism of \textsl{unbounded depth} it becomes possible to recruit arbitrarily many processors which can act in parallel to efficiently solve the sorts of ``combinatorily explosive'' problems warned about by Lighthill.   For this reason, the formal time metric for $\mathfrak{P}$ collapses \textsl{polynomial hierarchy} of complexity classes between $\mathbf{P}$ and $\mathbf{PSPACE}$.\footnote{See \citep{Savitch1979}.}   This contributes to a bulwark of arguments which suggest that concrete realizations of such ``unreasonable'' models cannot be constructed in practice so as to allow us to uniformly decide problems which cannot be decided efficiently using current computational architectures.\footnote{For instance \citep[p. 14]{Boas1990} states that the suggestion that we can build and apply devices allowing for unbounded parallelism entails  ``severe violations of basic laws of nature'' and that ``if physical constraints are taken into account, all gains of parallelism seem to be lost''.} 

All of this said, the current climate encourages speculation about whether such traditional barriers to feasible computation might be overcome by developments in artificial intelligence itself.  Popular accounts often given the impression that artificial intelligence systems already exploit parallelism which allows them to outperform ``classical'' techniques in various domains.  So one might wonder whether this will also turn out to be the case true for proof discovery.

The sort of machine learning systems currently under investigation -- e.g. recurrent networks or transformers -- do indeed exploit  ``parallel architectures'' in the sense that their operation typically relies on matrix operations which can be parallelized in principle.   In this sense, their high-level presentations are akin to programs for Parallel RAM machines.  On the other hand, using current or foreseeable computing machinery we are only able to directly carry out parallelized algorithms for (e.g.) matrix multiplication for inputs of a finite dimension which is fixed in advance by a concrete processor design.  As such, single applications of the systems in question cannot lead to \textsl{asymptotic} improvements of the sort of systems in question for arbitrarily large classes of instances as arise in the case of the proof discovery problem.\footnote{This remains the case even when the relevant sorts of techniques are carried out using specialized hardware such as \textsl{graphic} and \textsl{tensor processing units}.  These are optimized to perform matrix operation on low precision floating point numbers employed in machine learning applications but are otherwise of classical design.} As we will discuss further in \S 4.3, there thus appear to be principled reasons why machine learning techniques cannot on their own reduce the complexity of the discovery problem in the cases we have been considering.  

\section{Case studies}

To reiterate, the crux of the basic argument is that the proof discovery problem in mathematics is (provably) \textsl{computationally complex} and thus  \textsl{inherently difficult}.  But again, this argument concerns proof discovery understood as a decision problem with \textsl{infinitely many instances}.   On this understanding, it seems there is little room to challenge the argument on its own terms.   But we would be drawn to  question its relevance if the \textsl{specific cases} we cared about in practice ended up being significantly easier than the worst case reported by results like Theorems 1, 2, and 3.   It would, for instance, be highly germane if automated techniques already had  resolved a longstanding open question or appeared to be on the cusp of doing so.

At the time of writing, we are unaware of any instances in which the application of computing technology can be reasonably credited with resolving a \textsl{high profile} open question of the sort described in \S 2.  But there been several instances of lower profile results obtained in part or whole by the use of artificial intelligence-related methods.\footnote{\citep{Avigad2022} and \S 5 below takes step towards a systematic account of this distinction.} The best known examples to date have been obtained via traditional automated theorem proving techniques.   But there has also been a recent success employing $\textsc{Sat}$-solvers and large language models in a manner reminiscent of applications which are currently fueling public debates.

We will now examine three case studies illustrating the methods in question.  Each raises under-explored questions about the contours of mathematical knowledge and proof which may be of interest to philosophers in their own right.  But while illustrative, they are also sampled from a large body of evolving work whose details cannot be adequately surveyed here.  We will thus focus narrowly on demonstrating that these case studies exemplify the following characteristics:
\begin{enumerate}[(C1)]
\item The artificial intelligence-related methods which they employ are refinements of the technique traditionally known as \textsl{brute force search}.
\item The open questions which have been resolved by applying them are of a logical form (typically $\Sigma^0_1$) which allows us to see in advance that they are amenable to brute force while also being of a simpler logical complexity than the high-profile problems discussed in \S 2.  
 \end{enumerate}

\subsection{Automated theorem proving and the Robbins Problem}

Artificial intelligence and automated theorem proving have been  connected since their beginnings.  For instance at the 1956 conference where the expression ``artificial intelligence'' was coined, the only functioning system presented was the \textsl{Logic Theorist} of  Newell, Shaw, and Simon \citeyearpar{Newell1957}.  This was a program for proving statements of propositional logic in a manner intended to mimic how humans approached logic problems.  But as this method was incomplete,  Newell, Shaw and Simon themselves described it as \textsl{heuristic} and contrasted it with \textit{algorithmic} methods which are complete but need to return an output in a feasible number of steps.  

This situation was further addressed by \citet{Wang1960} who developed a sound and complete method for propositional derivability but abandoned the goal of modeling human reasoning.  Wang's approach was based on an operationalization of Gentzen's sequent in which the rules are both cut free and invertible.  This makes it possible to systematically perform proof search by working backwards from the main connective of the statement to be proved. Wang's method set the stage for the formulation of the \textit{Davis-Putnam algorithm} which serves as the basis of contemporary $\textsc{Sat}$-solvers of the sort discussed in \S 4.2.  But this method is itself based on the yet more fundamental rule of \textsl{resolution} which also remains at the core of most  automated theorem proving systems for first-order logic.\footnote{An accessible overview of the early history of automated theorem proving is provided by \cite[\S 3]{MacKenzie2001}.  See also \citep[\S 1]{Biere2021}.} 

Recall first that a propositional formula $\varphi_0$ can be efficiently transformed into a logically equivalent formula $\varphi_1$ in \textsl{conjunctive normal form} -- i.e. so that $\varphi_1$ is a conjunction of \textsl{clauses} each comprised of a disjunction of \textsl{literals} which are themselves either atoms $p$ or negated atoms $\neg p$.  If $\varphi_1$ contains distinct clauses of the respective forms $\psi \vee p$ and $\neg p \vee \chi$, then any valuation which makes $\varphi_1$ true must also make $\psi \vee \chi$ true.  This allows us to obtain a simpler formula $\varphi_2$ which is still equivalent to $\varphi$ in terms of satisfiability by removing all occurrences of $p$ from $\varphi_1$.  This idea is captured by the \textsl{resolution rule}:\begin{example}
\label{unitclause}
\begin{prooftree}
\AxiomC{$\varphi \vee p$}
\AxiomC{$\neg p \vee \psi$}
\RightLabel{\scriptsize $\mathscr{R}$}
\BinaryInfC{$\varphi \vee \psi$}
\end{prooftree}
\end{example}

As resolution-based theorem provers operate solely via this rule, they are simpler to implement than methods which have rules involving multiple propositional connectives  (like that of Wang).\footnote{In fact resolution is simply a form of the cut rule in the case the cut formula is an atom while Wang's method starts out with a cut-free formulation of the sequent calculus. See \citep[\S7.4]{Troelstra2000} for more on this relationship.} 
For this reason \citet[p. 24]{Robinson1965d} originally described resolution as  ``machine-oriented, rather than human-oriented''.   And it is indeed an application of this method which still provides one of the best-known applications of automated techniques to mathematics -- i.e. McCune's \citeyearpar{McCune1997} resolution of the \textsl{Robbins problem}.\footnote{Since this time, a number of other  results have been obtained in a similar manner which could also be used to illustrate the points framed below -- e.g. the recent resolution in a statement of non-commutative algebra known as the Weak Abelian Inner Mapping conjecture by Kinyon, Veroff, and Vojt{\v{e}}chovsk\`y. \citep{Kinyon2013} for a preliminary report).}

Recall that the theory $\mathsf{B}$ of Boolean algebras is given in the language $\mathcal{L}_{\mathsf{B}} = \{\wedge, \vee, \neg, 0, 1\}$ via universal closure of the following axioms:

\begin{tabular}{llcll}
(B$_1$)&$x \vee (y \vee z) = (x \vee y) \vee z$&~~~~&(B$_6$)&$x \wedge (y \wedge z) = (x \wedge y) \wedge z$\\
(B$_2$)&$x \vee y = y \vee x$&~~~~&(B$_7$)&$x \wedge y = y \wedge x$\\
(B$_3$)&$x \vee (x \wedge y) = x$&~~~~&(B$_8$)&$x \wedge (x \vee z) = x$\\
(B$_4$)&$x \wedge (y \vee z) = (x \wedge y) \vee (x \wedge z)$&~~~~&(B$_9$)&$x \vee (y \wedge z) = (x \vee y) \wedge (x \vee z)$\\
(B$_5$)&$x \vee \neg x = 1$&~~~~&(B$_{10}$)&$x \wedge \neg x = 0$
\end{tabular}
\vspace{0.1cm}

\noindent Inspired by an alternative axiomatization given by \cite{Huntington1933}, Robbins considered the following axioms $\mathsf{R}$ over the smaller language $\mathcal{L}^-_{\mathsf{B}} = \{\wedge,\vee\}$:
  
\begin{tabular}{llcll}
(R$_1$)&$x \vee (y \vee z) = (x \vee y) \vee z$&~~~~&(Def$_1$)&$0 = \neg (x \vee \neg x)$\\
(R$_2$)&$x \vee y = y \vee x$&~~~~&(Def$_2$)&$1 = x \vee \neg x$\\
(R$_3$)&$\neg(\neg(x \vee y) \vee \neg(x \vee \neg y)) = x$&~~~~&(Def$_3$)&$x \wedge y = \neg(\neg x \vee \neg y)$
\end{tabular}

\noindent 
It is not hard to see that R$_1$, R$_2$ and R$_3$ are derivable from $\mathsf{B}$.  Robbins posed the following problem: \textsl{if 0,1, and $\wedge$ are defined by \textnormal{Def}$_{1}$, \textnormal{Def}$_{2}$ and \textnormal{Def}$_{3}$, are all structures satisfying $\mathsf{R}$ Boolean algebras}?

Robbins and Huntington were unable to answer this question at the time and it was subsequently popularized by Tarski.\footnote{E.g., \citep[p. 245]{Henkin1971}.} Substantial progress would not be made until the late 1970s when Robbins' conjecture became the focus of the automated theorem proving group at the Argonne National Laboratory.  It was in this context that \citet{Winker1992} gave a conventional non-computer assisted proof that the adjunction of the single axiom $\omega := \exists{x}\exists{y}(x \vee y = x)$ to $\mathsf{R}$ allows for the derivation of each of B$_1$-B$_{10}$.  Resolving the Robbins problem was thus reduced to confirming
\begin{example}
$\mathsf{R} \vdash_{\mathrm{FOL}} \omega$
\label{winker}
\end{example}

This is what \citet{McCune1997} accomplished by employing a variant of the automated theorem OTTER called EQP.  Both of these systems take advantage of the fact that the resolution rule can be extended to the first-order case by combining it with a procedure for substituting terms for free variables to accomplish the task known as \textit{unificaiton}.\footnote{A 
\textsl{unifier} for expressions $E_{1}\sigma,E_{2}\sigma$ is a substitution $\sigma$ which makes them syntactically identical -- i.e. $E_{1}\sigma \equiv E_{2}\sigma$.  For instance if $E_1$ is $x = a$ and $E_2$ is $y = a$ then $\rho = \{b/x, b/y\}$ is a unifier since $(x = a)\rho \equiv b = a \equiv (y = a)\rho$ while $\sigma = \{y/x\}$ is a \textsl{most general unifier} as $\rho$ can be obtained from $\sigma$ by composing the latter with another substitution such as $\theta = \{b/y\}$. See, e.g., \citep[\S 7.2]{Troelstra2000} for more on the role of unification in resolution-based proof systems.}  OTTER also employs a general method for operationalizing equational reasoning known as the 
\textsl{Knuth-Bendix algorithm}.  This allows for the transformation of a given set of equations $\mathscr{E}$ into a confluent and terminating set of rewriting rules sufficient to check the truth of every equation in $\mathscr{E}$.\footnote{The basic idea is to observe that an equation between terms $t = s$ can be oriented so as to obtain a rewriting rule of the form $t \longrightarrow s$ (or $s \longrightarrow t$), allowing one to substitute the term on the left hand-side of the arrow with the term on the right hand-side. Given such a set of rules, a term $t$ may be repeatedly reduced to potentially obtain a non-reducible term $s$ which may then be equated with $t$ in further reasoning.   See, e.g., \citep{Dershowitz1990}.}  The method of resolution can further adapted in this context to treat equational reasoning exemplified by derivations from the axioms $\mathsf{R}$ via the technique of \textit{paramodulation}.\footnote{This takes the form
\begin{prooftree}
\AxiomC{$r = s \vee \psi$}
\AxiomC{$p[t] \vee \varphi$}
\BinaryInfC{$p\theta[s\theta] \vee \varphi \theta \vee \psi\theta$}
\end{prooftree}
where $\varphi$ and $\psi$ do not have any variable in common, $\theta$ is the most general unifier for $t$ and $r$, $p[t]$ indicates that the term $t$ appears in the atomic formula $p$, and $p\theta[s\theta]$ denotes the formula obtained by applying the substitution $\theta$ to both the atomic formula $p$ and the term $s$.   Paramodulation is thus a technique which combines unification, substitution, and resolution in a single rule.  See, e.g. \citep[\S 9.5]{Russell2021}.  The system EQP employed by McCune was further optimized to take into account the associativity and commutativity of the operation $\vee$ (as given by R$_1,R_2$).}    

McCune used this system to find two terms witnessing the existential quantifiers in (\ref{winker}) together with a proof witnessing the corresponding substitution instance.   These were ultimately quite short: $t_{1} := (x \vee x)$ and $t_{2} := \neg(\neg(x \vee x \vee x) \vee x)$) and a 13 line derivation of $t_1 \vee t_2 = t_1$ in the equational fragment of first-order logic.\footnote{The materials are all still available at \url{https://www.cs.unm.edu/~mccune/papers/robbins/}.}  But as we will discuss below,  the process of finding these objects was arduous in terms of computing resources.

This was heralded at the time as the first instance in which a computer had played a role in resolving a recognized open question in a manner which exhibited ``creativity'' or ``logical insight''.\footnote{E.g. \citep{Kolata1996}.}  Nonetheless, the manner in which McCune obtained his result exhibits two characteristics which are shared by several subsequent examples in automated theorem proving.  First, the statement which was derived by EQP was not a direct statement of the original problem -- i.e. \textsl{all are Robbins algebras Boolean?} -- but rather a reformulation in terms of formal derivability.  Second, the sufficiency of this formulation for resolving the problem was not immediate given its original statement but rather depended on an intervening lemma which was obtained by a human mathematician.

The significance of these characteristics comes into focus when we further examine the details of how EQP operates.  Since this system descends from OTTER, it implements a sound and complete proof procedure for first-order logic.  However its main loop simply dovetails  the derivation of formulas from axioms (via resolution), the inductive generation and simplification of $\mathcal{L}_{\mathsf{B}}$-terms (via the Knuth-Bendix algorithm), and the simplification of derived formulas (via paramodulation).   

Per claim (C1), it is  his latter observation which illustrates  what we mean by saying that the method employed in McCune's proof is an instance of \textsl{brute force search}.  Such a procedure is traditionally understood as a means of solving a decision problem $X$ in cases where positive instances $x \in X$ are witnessed by the existence of so-called \textsl{certificates} $c$ drawn from an effectively enumerable class $C_X$ satisfying a decidable predicate $\Phi(x,c)$. Positive instances of problems which admit such a characterization are thus verifiable by a so-called \textsl{generate-and-test} procedure which for a given instance $x$ simply enumerates the members of $C_X$ as $c_0,c_1,c_2,\ldots$ and checks whether $\Phi(x,c_i)$ holds.\footnote{See \citep{Trakhtenbrot1984}, \citep{Heule2017}, and \citep{Dean2019c} for more on the history and significance of the notion of brute force in computer science.}    

Standard techniques allow us to see that any decision problem $X$ admitting such a characterization can be defined by a $\Sigma^0_1$-formula in the language of first-order arithmetic -- i.e. so that $X = \{x : \exists m \varphi(\ulcorner x \urcorner,  m) \}$ for $\varphi(x,y)$ a $\Delta^0_0$-formula and $\ulcorner \cdot \urcorner$ an appropriate arithmetical encoding.  Per claim (C2) this shows that the logical form of the Robbins problem may naturally be regarded as $\Sigma^0_1$-statement of  $\mathcal{L}_a$ -- a feature which we will suggest in \S 5 classifies it as \textsl{simpler} than the sorts of high-profile open questions considered in \S 2.\footnote{Of course questions about formal derivability from a computable set of axioms $\Gamma$ -- of which (\ref{winker}) is typical -- are paradigmatically of this sort as in this case we may take $\Phi(x,y)$ to be a formalized proof predicate $\textrm{Proof}_\Gamma(x,y)$ for $\Gamma$ and the relevant class of certificates $C_{\Gamma}$ is the set of well-formed derivations from $\Gamma$.  If we equate proof discovery with formal derivability in the manner we have proposed in \S 2, then the solution to \textsl{any} open question can be reduced to an enumerative search from axioms as just described.  We will return to this issue in \S 5.}

In addition to this, McCune's proof also possesses several other features which are traditionally associated with brute force.  First, note that the characterization just given allows in general that $C_X$ may be infinite.  As a consequence, the corresponding generate-and-test procedure for $X$ will not terminate if $x \not\in X$.  This in turn corresponds to the familiar observation that if $\Gamma \not\vdash_{\mathrm{FOL}} \varphi$, then an automated theorem prover will not be able to detect this by carrying out an enumerative search through formal derivations.  Thus had Robbins' conjecture been \textsl{false} -- in which case $\mathsf{R} \not\vdash\omega$ -- then a system like EQP would have been unable to resolve the problem in the manner of McCune's proof.%\footnote{The falsity of Robbins' conjecture is equivalent to the statement that there is a counterexample in form of a structure $\mathcal{A}$ satisfying the axioms $\mathsf{R}$ but not $\mathsf{B}$.   There is no in principle obstacle to formalizing this statement itself as a formula $\psi$ which could in principle itself be derived by the an automated theorem prover.   Note, however, a formal language in which $\psi$ could be expressed would have to be more expressive than that of $\mathcal{L}_{\mathsf{B}}$. In particular, since $\mathcal{A}$ might be infinite, $\psi$ would most naturally be taken as either a $\Sigma^1_1$-statement in the language of second-order arithmetic or a similarly complex statement in the language of first-order set theory.   But proving such statements via automated techniques in general requires methods which go considerably beyond those implemented in EQP (or related finite model checking systems).} 

Second, in addition to our prior observation about the complexity of the meta-theoretic assertion  (\ref{winker}), the \textsl{formula} $\omega$ is itself a $\Sigma_1$-formula in the language of Boolean algebras.   As such, one might wonder whether $\omega$ can only be derived in $\mathsf{R}$ is via non-constructive proof which makes essential use of quantificational reasoning.  But as we have seen, McCune's proof did not proceed in this manner.   Rather it employed an inductive generation procedure within EQP to enumeratively  search for both $\mathcal{L}^-_{\mathsf{B}}$-terms and a derivation witnessing $\mathsf{R} \vdash t_1 \vee t_2 = t_1$.

The manner of this search was not entirely ``naive''. But third, it may still be reasonably characterized as \textsl{general} or \textsl{domain unspecific}.   This is to say that the procedures in question were in principle capable of generating and substituting any well-formed terms and deriving any equational proof.  But the order in which these steps were carried out in practice was not informed by the \textsl{mathematical content} of the Robbins problem or facts about Boolean algebras in general but rather by general-purpose parameters controlling the operation of EQP.\footnote{For instance McCune describes how the software was configured to vary both the length of terms which were generated and also how terms of different ``ages'' were substituted into previously generated equations.} 

Fourth, while a modest sized proof and terms were eventually found, the overall enumerative procedure by which they were obtained can reasonably be described as proceeding through an \textsl{unsurveyably large} number of intermediate steps.\footnote{McCune reports that during the execution which ultimately yielded the result, EQP generated 49,548 equalities from 2,612,977 attempts at term rewriting which in turn required 8 days of computing time and used about 30 megabytes of memory.   But this was only obtained after a multi-week process of failed executions with different parameter values of the sort described in the prior note.}  In other words, while we can comprehend the terms and proof described above after the fact, it would have been far beyond the practical abilities of a human mathematician working without automation to find them in the manner of McCune's proof.

\subsection{\textsc{SAT} solvers and finite combinatorics}
\label{sec4.2}

Recall that $\textsc{Sat}$ denotes the satisfiability problem for propositional logic -- i.e. that of checking if a formula is true in \textsl{some} row of its truth table.  A $\textsc{Sat}$-\textsl{solver} is a decision algorithm satisfying:
\begin{example}
\label{satprops}
\begin{enumerate}[i)]
\item \textsl{Totality}: $\alpha(\varphi)$ returns  output $1$ -- i.e. ``satisfiable'' -- or $0$ -- i.e. ``unsatisfiable'' -- for all propositional formulas $\varphi$.
\item \textsl{Soundness}: If $\alpha(\varphi) = 1$, then $\varphi \in \textsc{Sat}$.
\item \textsl{Completeness}: If $\alpha(\varphi) = 0$, then $\varphi \not\in \textsc{Sat}$ -- i.e. $\varphi \in \overline{\textsc{Sat}}$.
\end{enumerate}
\label{tsc}
\end{example}
Recall also that if $\varphi$ is \textsl{unsatisfiable} -- i.e. \textsl{false} in all rows of its truth table -- then $\neg\varphi$ is a tautology and thus provable from no premises in the propositional calculus (by the completeness theorem). As such, an algorithm $\alpha$ satisfying (\ref{satprops}i-iii) also serves as a decision method for  $\textsc{Prov}_{\emptyset}^{\mathrm{PL}}$.   

As we have discussed, $\textsc{Sat}$ is $\mathbf{NP}$-complete.  Presuming that $\mathbf{P} \neq \mathbf{NP}$, there thus cannot exist an algorithm satisfying (\ref{satprops}i-iii)  implementable by a reasonable model of computation which always returns an output in time polynomial in the number of propositional variables in $\varphi$.  And in fact no procedure with better than exponential worst case running time in the general case is known to exist. It is thus striking that a number of algorithms which perform well on large classes of $\textsc{Sat}$ instances are currently being investigated.\footnote{Although $\mathbf{P} \neq \mathbf{NP}$ entails that there cannot exist a polynomial time algorithm for $\textsc{Sat}$, it is perhaps already surprising that there are decision algorithms for $\textsc{Sat}$ which are asymptotically more efficient than the naive $O(2^n)$ method of truth tables (see \citealp[\S 16]{Biere2021}).}      

Most contemporary $\textsc{Sat}$ solvers are refinements of the so-called Davis-Putnam algorithm.  This is itself a recursive implementation of the method of resolution described in \S4.1 together with a strategy for selecting a literal on which to apply a rule known as \textsl{unit propagation}.\footnote{For instance, consider the conjunctive normal form formula $\varphi_0 := (p \vee q) \wedge (\neg p \vee r) \wedge (\neg r \vee s) \wedge p$.  This contains a so-called \textsl{unit clause} $p$ consisting of a single literal.   As a valuation making $\varphi_1$ true must make $p$ true, $\varphi_0$ can be simplified by removing any clause contain $p$ and deleting $\neg p$ from any clause in which it occurs to yield a simplified formula $\varphi_1 := r \wedge (\neg r \vee s)$ which is equivalent in terms of satisfiability.  This process can then be re-applied to the unit clause $r$ to obtain $\varphi_2 := s$}  However this is not a fully explicit algorithm.   Specific $\textsc{Sat}$ solvers thus rely on a combination of heuristics for choosing literals according to the properties of the input formula or backtracking when conflicts are found.  Appropriate choices of this sort lead to algorithms which reduce the overall number of assignments which must be considered in large classes of cases.\footnote{See, e.g., \citep[\S 7]{Russell2021} for an overview and \citep{Biere2021} for a presentation of the specific techniques used in the Cube and Conqueror algorithm described below.}

The applicability of $\textsc{Sat}$ solvers to  proof discovery arises due to the possibility of using propositional logic to express mathematical statements either directly or as a parameterized family which is successively checked by an enumerative search.   The process can be illustrated by examining one of the successful applications of such algorithms to answer an open question -- i.e. the recent proof of a statement in  combinatorics known as the Boolean Pythagorean Triple Conjecture proposed by Ronald Graham in the 1980s for which he offered a \$100 prize.\footnote{BPT is similar in form to earlier results in additive combinatorics such as Schur's Theorem on monochromatic sums (1912), van der Waerden's Theorem on monochromatic arithmetical progressions (1927), and Szemer\'edi's Theorem (1975) on arithmetical progressions in initial segments of positive upper Banach density.  See \citep{Arana2015} for account of the \textsl{mathematical depth} of such related this is related to our disccussion in \S 5.}  

A \textsl{Pythagorean triple} is $a,b,c \in \mathbb{N}^+ = \{1,2,3,\ldots\}$ such that $a^2 + b^2 = c^2$ -- e.g. $3,4,5$ and $5,12,13$.   Graham's conjecture was that for every $k \geq 1$, there exists an $n$ such that for all $k$-colorings $f(x)$ of $\{1,\ldots,n\} =_{\mathrm{df}}[n]$ -- i.e. mapping $f:[n] \rightarrow \{0,\ldots,k-1\}$ -- contains a so-called \textsl{monochromatic} Pythagorean triple -- i.e. $a,b,c$ which are assigned the same value by $f(x)$.  The Boolean Pythagorean Triple Conjecture is the case of this question for $k = 2$ -- i.e.
\begin{itemize}
\item[(BTP)] $\exists n \in \mathbb{N} \ \forall f:[n] \rightarrow \{0,1\} \ \exists a,b,c \in [n] \textrm{ such that } a^2 + b^2 = c^2 \ \& \ f(a) = f(b) = f(c)$
\label{bpt}
\end{itemize}
If BPT is true,  there is a least witness $m \in \mathbb{N}$ such that
\begin{itemize}
\item[(BTP$_m$)] $\forall f:[m] \rightarrow \{0,1\} \ \exists a,b,c \in [n] \textrm{ such that } a^2 + b^2 = c^2 \ \& \ f(a) = f(b) = f(c)$
\label{bpt}
\end{itemize}

The recent work of \citet{Heule2016} completely resolved BPT by showing that that while BPT$_{7824}$ is \textsl{false}, BPT$_{7825}$ -- and hence also BPT itself -- is \textsl{true}.\footnote{This is one of several results which have recently been been obtained using $\textsc{Sat}$-solvers -- e.g. the resolutions of the Keller conjecture by \citep{Brakensiek2020}, the calculation of the packing number of the infinite square by \citep{Subercaseaux2023} and the number of points in general position required to ensure the existence of a convex hexagon without an interior point (an instance of Erd\H{o}s's \textsl{happy endings problem}) by \citep{Heule2024}.   These cases use more sophisticated propositional encodings which illustrates the division of human and automated labor  which is often required to apply $\textsc{Sat}$-solving techniques to the original formulation of combinatorial problems.  But they could otherwise be used to frame similar points.}  Both of these facts were demonstrated by using a $\textsc{Sat}$ solver known as \textsl{Cube and Conquer} [C\&C].  The first fact was shown by using this algorithm to find a $2$-coloring $f(x)$ of $\{1,\ldots, 7824\}$ such that neither of the sets respectively colored $0$ or $1$ by $f(x)$ contains a Pythagorean triple. The second fact was shown by using C\&C to demonstrate that no such coloring of $\{1,\ldots, 7825\}$ exists.   

Since BPT is a statement involving natural numbers and the inputs to C\&C are propositional formulas, it might at first seem mysterious how this was possible.  But it is straightforward to see how to construct for any fixed $m$ a propositional formula such that  BPT$_m$ is true if and only if $\beta_m$ is \textsl{not} satisfiable.\footnote{For let $T_m = \{\langle a,b,c \rangle \in [m] : a^2 + b^2 = c^2\}$ be the set of Pythagorean triples up to $m$ and consider propositional formulas of the form $\tau_{a,b,c} = (x_a \ \vee \ x_b \ \vee \ x_c) \ \wedge \ (\neg x_a \ \vee \ \neg x_b \ \vee \ \neg x_c)$ for $\langle a,b,c \rangle \in T_m$.  Note also that any 2-coloring $f:[m] \rightarrow \{0,1\}$ can also be viewed as a propositional valuation $v_f(x_i)$ by interpreting the truth/falsity of the propositional variable $x_i$ as expressing that the number $i$ is assigned the color $1/0$ by $f(x)$.   It follows that $v_f(\tau_m) = 1$ if and only if $f$ is a coloring for which $f(a), f(b), f(c)$ are not all assigned the same value.  It thus also follows that
\begin{example}
$\beta_m = \bigwedge_{\langle a,b,c \rangle \in T_m} \tau_{a,b,c}$ 
\end{example}
is true relative to $v_f$ if and only if $f$ is a $2$-coloring of $[m]$ which does \textsl{not} contain a monochromatic Pythagorean triple.  Hence $\beta_m$ is \textsl{unsatisfiable} just in case BPT$_m$ is true.}
The number of propositional variables in $\beta_m$ is given by the number of distinct natural numbers which appear in Pythagorean triples with members $\leq m$.  After taking account  easily identified symmetries, it is possible to construct formulas $\beta_{7824}$ and $\beta_{7825}$ which respectively contained 3730 and 3745 variables.     

Testing formulas of this size for satisfiability by truth tables is indeed highly infeasible using even remotely foreseeable computing technology.  It is thus testament to the advances in $\textsc{Sat}$ solving techniques that \citet{Heule2016} were able to concretely apply the C\&C algorithm to these formulas to demonstrate that $\beta_{7824} \in \textsc{Sat}$ -- and thus BTP$_{7824}$ is \textsl{false} -- and that $\beta_{7825} \not\in \textsc{Sat}$ -- and thus BTP$_{7825}$ and also BTP itself are \textsl{true}.   This yielded a so-called \textsl{certificate of unsatisfiability} which, while not a conventional deductive proof itself, can be efficiently converted into a formal refutation of $\beta_{7825}$ in (e.g.) a traditional Hilbert system.   On the other hand, these objects are very large in the sense of requiring much memory (over 200 terabytes).   Finding them also required much computational labor to obtain (over 4 CPU years) using a powerful computer (800 cores).  

In this sense, the proof of the Boolean Pythagorean Triple Conjecture is akin to other ``computer proofs'' which required extensive computational effort as well as resulting in large proof objects.   Perhaps the best known of these is Appel and Haken's \citeyearpar{Appel1977a} proof of the Four Color Theorem in which a computer was originally employed to exhaustively confirm a large case distinction whose verification required confirming that over 4000 finite graphs satisfied a decidable property.    But in this case the role of the computer was confined to confirming a lemma in an otherwise human-generated proof strategy.  

On the other hand, almost all of  Heule et al.'s \citeyearpar{Heule2016} proof was computer generated.  In particular, its final step consisted in applying the general purpose C\&C  algorithm to the formula $\beta_{7825}$.  Nonetheless, they also formally verified their entire proof of BPT -- inclusive of the correctness of C\&C  -- in the manner which has now been applied to the  Four Color Problem and the Kepler Conjecture.\footnote{See \citep{Gonthier2013} and \citep{Hales2017}.}  There can thus be no question as to whether Heule et al.'s method in fact discovered a \textsl{correct} proof of BPT .   Nonetheless, the demonstration in question was again far from being surveyable by a human mathematician.\footnote{See \citep{Tymoczko1979} for the classical discussion of the concept of surveyability in the philosophy of mathematics.   The analysis is incisively refined by \citep{Heule2017} who, amongst other things, conjecture that BPT is what they refer to as an ``alien truth'' -- i.e. one whose statement is sufficiently short so that it can be grasped by humans but for which it seems likely that there exist only unsurveyably large brute-force like derivations.}   

Further to (C1), \citet{Heule2017} explicitly observe that this is in large part due to the fact that $\textsc{Sat}$-solvers like C\&C are again implementations of brute force search.   This is to say that in order to decide whether a particular formula $\beta_m$ is satisfiable, they exhaustively explore different ways in which it can be made true.  This process was implemented in a sophisticated way so that in some cases parts of the search space may be excluded -- e.g. in virtue of the distribution of propositional atoms across different clauses in the input formula.  But like traditional automated theorem provers,  the method of C\&C is general in the sense that it does not take into account what we might otherwise call the \textsl{mathematical content} expressed by the formulas to which it is applied.

In addition to this, the  method of \citep{Heule2016} also exhibits another characteristic of brute force.   For since it was initially unknown whether there existed $m \in \mathbb{N}$ for which $\mathrm{BPT}_m$ holds, their proof required the successive application of C\&C to a family of formulas $\beta_i, \beta_{i+1}, \beta_{i+2}, \ldots$ The overall structure of their proof thus implemented an unbounded search whose termination could not be guaranteed \textsl{a priori}.   Thus had BTP (or another similarly encoded combinatorial statement) been false, it could not have been refuted in this manner.   

The status of (C2) in regard to BTP is in turn related not just to its logical form, but also that of the class of formulas which can be decided by $\textsc{Sat}$-solvers in general.  Note first that it is again easy to see from (\ref{bpt}) that $\mathrm{BTP}$ is equivalent to a $\Sigma^0_1$-statement of $\mathcal{L}_a$ --  i.e. one of the form $\exists x \varphi(x)$ with a unbounded existential quantifier and a decidable matrix $\varphi(x)$.  In fact, this definition shows that it is equivalent to such a formula in special case where $\varphi(x)$ can be written as a so-called $\Pi^b_1$-\textsl{predicate} --  i.e. so that $\varphi(x)$ has the form $\forall x < t \psi(x)$ where $t$ is a term in the \textsl{language of bounded  arithmetic}.    Since the sets of natural number definable by $\Pi^b_1$-predicates correspond to the problems in $\mathbf{coNP}$, this gives a precise syntactic characterization of the class of formulas which are decidable in general by a single application of a $\textsc{Sat}$-solver.\footnote{See \citep{Buss1998b} for the definition of $\Pi^b_1$ as well as the proof that formulas define the precisely the languages in the class $\mathbf{coNP}$.}

\subsection{Large language models and the cap set problem}

It is a commonplace that a distinction should be drawn between automated theorem proving and machine learning techniques.   But such a classification falls short of a precise characterization of the difference between what we will call \textsl{logic-based} and \textsl{statistics-based} methods.  This in turn complicates systematically addressing the topical question: \textsl{How might the latter class novelly contribute to proof discovery beyond the former?}

In approaching this we may also distinguish between what might be called \textsl{systematic} and \textsl{non-systematic} applications of statistic-based methods to proof discovery.   The former seek to find applications of machine learning within the framework of automated theorem proving surveyed in \S 4.1. This is exemplified by attempts to refine the so-called method of \textsl{hammers} -- i.e. general purpose techniques to fill in gaps in proofs by consulting libraries of previously formalized theorems in a manner reminiscent of expert systems (see, e.g., \citealp{Blanchette2016}).   

Such libraries are now sufficiently large to support an ongoing project to improve upon the performance of hammers by applying machine learning to the so-called \textsl{premise-selection problem} -- i.e. that of proposing candidate lemmas to prove a give target theorem based on the statistical properties of the linguistic structure of the theorem they contain within a corpus of previously formalized proofs.    Techniques of this sort have become increasingly successful in improving the performance of automated theorem provers in rediscovering fully automated proofs of previously verified theorems.\footnote{See, e.g., \citep{Urban2013} and \citep{Hales2014}.}  But at the time of writing (and to the best of our knowledge) this method has not been used to resolve any standing open questions.

This stands in contrast to the use of machine learning which has recently figured in what is the first successful application of a paradigmatic statistics-based method to proof discovery (again to the best of our knowledge).\footnote{Needless to say, other related attempts are currently underway -- see, e.g., \citep{Lin2024}. Other examples may well be available by the time this paper appears.}  This comes in the form of using a large language model to aid in the discovery of a structure known as a \textsl{cap set} which improved upon known lower bound results.   We will again summarize the mathematical result before describing the methods which were recently used by \citet{Romera-Paredes2024} to obtain it in regard to (C1) and (C2).     

The notion of a cap set also originates in arithmetical combinatorics.    Recall that $\mathbb{Z}/m$ denotes the integers $\mathrm{mod}\ m$ considered as a cyclic group.  $(\mathbb{Z}/m)^n$ in turn denotes the set of vectors of length $n$ formed of elements of $\mathbb{Z}/m$.  A \textsl{cap set} is a subset of distinct elements of $(\mathbb{Z}/m)^n$ which  do not form an arithmetical progression.\footnote{For instance, in the case $m = 3$ and $n = 4$ a cap set $C \subseteq \{0,1,2\}^4$ consisting of vectors of the form $\vec{v} = \langle v_0,v_1,v_2,v_3 \rangle$ such that there do not exist three distinct $\vec{x},\vec{y},\vec{z} \in C$ such that $\vec{y} - \vec{x} = c(\vec{z} - \vec{y})$ for $c \in \{0,1,2\}$.   This condition can be seen to be equivalent to $\vec{x} + \vec{y} + \vec{z} \neq 0$ where the addition is performed $\mathrm{mod}\ 3$.  See, e.g., \citep{McMahon2019} for an acesssible introduction to this topic in relation to the card game SET.} Background developments motivate the study of the function $\mathrm{cap}(n) =$ \textsl{the largest subset of $(\mathbb{Z}/3)^n$ without 3-term arithmetic progressions}.  While exact values of this function are known only for $n \leq 6$, it is straightforward to see that $2^n$ provides a lower bound on $\mathrm{cap}(n)$ and that $3^n$ is an upper bound.  The original \textsl{Cap Set Conjecture} was as follows:
\begin{example} Does there exist $c < 3$ such that $\mathrm{cap}(n) \leq c^n$?
\end{example}

%2^3×3^2×11×71

This question was open for 20 years and was regarded as both difficult and significant within number theory at large.\footnote{E.g. because of its relation to Erd\H{o}s and Tur\'an's longstanding conjecture that if $A \subseteq \mathbb{N}$ is s.t. $\Sigma_{n \in A}\frac{1}{n} = \infty$, then $A$ contains arbitrarily long arithmetical progressions.  See, e.g., \citep{Grochow2019}.}  It was resolved in the positive by \citet{Ellenberg2017} who showed $\mathrm{cap}(n) \leq 2.756^n$ via a traditional non-computer assisted proof.  The longstanding lower bound of $2.2173^n \leq \mathrm{cap}(n)$ was due to \citet{Edel2004} who described a uniform construction of cap sets which operates by recursively composing known cap sets in smaller values of $n$ with so-called \textsl{admissible sets} so as to obtain cap sets for larger values.   A consequence of Edel's method is that discovery of larger admissible sets -- which are themselves subsets of $(\mathbb{Z}/3)^n$ -- can be used to construct larger cap sets at arbitrarily large higher dimensions.  Edel's lower bound was recently improved to $2.2180^n$ by \citet{Tyrrell2023} who employed a $\textsc{Sat}$-solver to search for such admissible sets.\footnote{See \citep{Tyrrell2023} Definition 2.8 for the definition of admissibility and lemmas leading to the proof his Theorem 1.2 on which the discussion below is based.}

The foregoing reprises the situation prior to the yet more recent work of \citet{Romera-Paredes2024}.  Their contributions can be summarized as follows:
\begin{example}
\label{cap}
\begin{enumerate}[i)]
\item $\mathrm{cap}(8) \geq 512$ 
\item $2.2194^n \leq \mathrm{cap}(n)$
\end{enumerate}
\end{example}  
(\ref{cap}i) improves upon the prior lower bound $\mathrm{cap}(8) \geq 496$ and and was established by finding a specific cap set $R$ of size 512 in $(\mathbb{Z}/3)^8$.  (\ref{cap}ii) concerns the asymptotic behavior of $\mathrm{cap}(n)$ for all $n$.  But it too was obtained by finding a single finite object -- an admissible set called $\mathcal{I}(15,10)$ -- which can be used to construct a family of cap sets in a higher dimension in a manner similar to Edel's construction.\footnote{\citep{Romera-Paredes2024} report a stronger result based on the construction of what they call a \textsl{partial admissible set} $\mathcal{A}(24,17)$. But the calculation on which this bound is based is not stated explicitly.  We will thus consider instead their construction of $\mathcal{I}(15,10)$ as it can be directly compared with that of admissible set $\mathcal{I}(11,7)$ given by \citep{Tyrrell2023}.} 

Before examining question (C1) in this case, let us first address the logical complexity of statement (\ref{cap}) into regard to (C2).  Since cap sets are finite  objects, it is easy to see that the property of being a cap is a decidable predicate.   The statement that $R \subseteq (\mathbb{Z}/3)^8$ is a cap set of size 512 is thus expressible as a $\Delta^0_1$-formula in the language of first-order arithmetic -- i.e. one provably equivalent to both a $\Sigma^0_1$- and $\Pi^0_1$-formula.  But also note that since $|(\mathbb{Z}/3)^n| = 3^n$, there are only finitely many  $X \subseteq (\mathbb{Z}/3)^n$ which are potential cap sets.  As such, the statement $\mathrm{cap}(8) > 496$ is equivalent to $\exists X \subseteq (\mathbb{Z}/3)^8(496 < |X| \leq  6561 \ \wedge \ X \text { is a cap set})$.  

Since finite sets can be coded as single natural numbers, this allows us to express (\ref{cap}i) itself as a $\Delta^0_1$-sentence of $\mathcal{L}_a$.  As stated, (\ref{cap}ii) is equivalent to a $\Pi^0_1$-sentence of $\mathcal{L}_a$.  But as we have already noted, it was obtained not by proving the universal statement $\forall n(2.2194^n \leq \mathrm{cap}(n))$ directly -- e.g. by  induction on $\mathbb{N}$ -- but rather by showing that $\mathcal{I}(15,10)$ is an admissible set. Thus to improve the prior lower bound on $\mathrm{cap}(n)$, a search over potential admissible set of increasing size was required.   Since this search was \textsl{a priori} unbounded, the statement which was in effect demonstrated by \citet{Romera-Paredes2024} is equivalent to a $\Sigma^0_1$-sentence of $\mathcal{L}_a$.   

These observations are in keeping with the means by which \citet{Romera-Paredes2024} demonstrated (\ref{cap}i,ii). This involved a combination of the techniques under the headings of \textsl{greedy algorithms}, \textsl{genetic algorithms}, and \textsl{large language models}.\footnote{See, e.g., \citep[\S 16]{Cormen2005} and \citep{Holland1992} for an account of the first two of these methods (which are both ``classical'' techniques not directly related to artificial intelligence in inception).}   In order to understand why these techniques needed to be applied in concert, note first that the problem of finding a cap (or admissible) set in $(\mathbb{Z}/n)^3$ can in principle be solved by searching through the $2^{(3^n)}$ subsets of $X \subseteq \{0,1,2\}^n$.  As this is (highly) infeasible already for $n = 8$ -- in which case more than $10^{1975}$ sets would need to be checked -- a more efficient strategy was evidently required.  

One means of attacking the problem is to search not for a cap set itself, but rather for a smaller description of such an object in the form of a procedure which outputs a potential cap set $X \subseteq \{0,1,2\}^n$ (from no input). One way in which such a procedure might operate is to follow the template of a greedy algorithm $\gamma$ -- i.e. one which attempts to construct $X$ in stages $X_0 \subset X_1 \subset X_2 \subset \ldots$ such that $X_{i+1}$ is obtained from $X_i$ by adjoining a new vector $\vec{v}$ just in case it does not form an arithmetical progression along with any two vectors already in $X_i$.\footnote{An antecedent to this general method is provided by the ongoing project of \citep{Gauthier2023a,Gauthier2023b} to use tree neural network to discover Python-like programs which compute previously distinguished entries in the \textsl{Online Encyclopedia of Integer Sequences} \citep{Sloane2007} via novel formulas.}    

To implement this template a \textsl{priority function} $p:(\mathbb{Z}/3)^n \rightarrow \mathbb{R}$ must also be supplied which returns a value ranking vectors determining the order in which they are considered for adjunction to $X_i$.  Such a function may  itself be specified by a \textsl{priority program} $\rho$ (in this case given concretely as Python code) which computes a function of the appropriate type.  We thereby obtain a program $\gamma[\rho]$ which may be executed to obtain a uniquely determined cap set $X_{\gamma[\rho]}$ in $(\mathbb{Z}/3)^n$.  

We can then define the \textsl{score} $s(\rho) = |X_{\gamma[\rho]}|$ -- i.e. the size of the cap set returned by $\gamma[\rho]$.  In the terminology of genetic algorithms, this is analogized to the \textsl{fitness} of $\rho$ in generating cap sets.  Of course an arbitrary $\rho$ is unlikely to yield an algorithm $\gamma[\rho]$ which produces a large cap set.  But genetic techniques also provide specific strategies for  attempting to improve the score of $\rho$ by treating it as a member of a \textsl{population} of high-scoring priority programs $\Pi$.  These were sampled by a \textsl{selection method} $\sigma$ which iteratively culls programs according to their fitness.   

In the next stage, these programs were encoded as \textsl{prompts} to a large language model.   Such models may in general be understood as a neural network derived by applying a learning algorithm $\lambda$ to a corpus of text $\mathcal{C}$ in a manner which results in an object $M$ (similar to a matrix of weights) such that when $M$ is presented with a prompt $p$ it produces an output $M(p) = q$.  On the intended interpretation, $q$ is be understood to be \textsl{similar} to the prompt $p$ in the sense of the statistical regularities encapsulated in $\mathcal{C}$ which have been distilled via $\lambda$ to produce $M$.  In the specific case in question, the prompt to $M$ was of the form $p = \langle \rho_1, \rho(\rho_1), \rho_2, s(\rho_2) \rangle$ consisting of the two high scoring programs produced at the prior step appended with their scores.   The output of $M(p)$ was next treated as a priority program with the hope that it might score higher than  $\rho_1$ and $\rho_2$.  If so, it was then added to the population $\Pi$. The prior steps were then iterated until a cap set in the given dimension of sufficiently large size was found.   

These steps comprise the method which \citet{Romera-Paredes2024} dub \textsl{FunSearch}.  It is evident that many highly specific decisions about $\gamma, s, \Pi,\sigma,\mathcal{C}, M$, and $\lambda$ are required to obtain a piece of executable software from the foregoing description.  This being said, \cite{Romera-Paredes2024} describe their result as follows: 
\begin{quote} \footnotesize{
Our proposed method, FunSearch, pushes the boundary of LLM-guided evolutionary procedures to a new level: the discovery of new scientific results for established open problems $\ldots$ Surpassing state-of-the-art results on established open problems provides a clear indication that the discoveries are truly new, as opposed to being retrieved from the LLM's training data. (p. 468)}
\end{quote}

Such claims immediately raise (at least) the following questions:

\noindent \textsl{Taking into account the disparate methods which contribute to FunSearch, what was the \textnormal{specific role of machine learning} in obtaining \textnormal{(\ref{cap}i,ii)}? Could these results have been demonstrated in another way using other extant or foreseeable techniques?  What are the prospects for applying the techniques in question to other open questions?}

\noindent A thorough assessment of these issues is beyond our current scope. What we will concentrate on is how they relate to (C1) -- i.e. in this case the claim that the relevant application of \textsl{FunSearch} represents another instance of brute force search.   

The fact that such a characterization is appropriate becomes apparent as soon as it is taken into account that the main loop of FunSearch had to be iterated \textsl{millions of times} to obtain the results (\ref{cap}i,ii).  In particular, by no means did the large language model $M$  immediately produce a priority program $\rho$ such that the resulting greedy algorithm $\gamma[\rho]$ generated the relevant cap or admissible sets.  Rather it was only as a consequence of successively using $M$ to suggest a long sequence of priority programs which were then updated according to a specific genetic algorithm that the appropriate cap and admissible sets were found.  

On this basis it is possible to record several similarities between the results of \citep{Romera-Paredes2024} and those surveyed in \S 4.1-\S 4.2.   First, we have already noted that that the statements (\ref{cap}i,ii) can both be regarded as $\Sigma^0_1$-statements which were demonstrated to be true by finding a concrete witness.   Second, the manner in which this witness was found in both cases was by a form of enumerative search -- i.e. potential witnesses were generated by the procedure just described and then tested to see if they had the desired property of being a sufficiently large cap or admissible set to surpass the previously known bounds.\footnote{Note, however, that there is a substantial difference in content between (\ref{cap}i,ii) and Heule et al.'s solution to the BPT Conjecture surveyed in \S 4.2.   In particular, while both results were obtained by an enumerate-and-test procedure, \citep{Heule2016} not only demonstrated that BTP was true, they determined that 7825 was the \textsl{precise lower bound} on the size of a set for which every two coloring contains a monochromatic Pythagorean triple.   This in turn required using the C\&C solver to demonstrate something which is provably equivalent to the fact that $\neg \beta_{7825}$ is a \textsl{tautology} -- i.e. true in \textsl{every row} of its truth table.  On the other hand, the results of \citep{Romera-Paredes2024} simply improved on known lower bounds which could in principle also have been obtained in the manner of Tyrrell's \citeyearpar{Tyrrell2023} application of a $\textsc{Sat}$-solver to the (presumptively) easier task of showing that an appropriately constructed propositional formulas is \textsl{satisfiable} -- i.e. true in \textsl{some} row of its truth table.  The authors do not address the question of whether FunSearch is capable of surveying the entire space of priority programs.   But as we discuss in note \ref{kolmogorovnote}, there is in fact reason to suspect that there are simple priority programs which such a method will not discover.} Third, the number of cases which had to be tested was unsurveyable from a human perspective.  Fourth, the overall search method was \textsl{general} or \textsl{domain unspecific}.

In regard to the last of these observations, two points stand out.  First, while \citet{Romera-Paredes2024} indicate that approximately 2.5 million iterations of the loop described above were required to obtain (\ref{cap}i,ii), this is far smaller than the number of combinatorial possible cap sets $X \subseteq \{0,1,2\}^8$.  In light of this it might initially seem appropriate to describe $M$ as having ``understood'' something about the mathematical content of the cap set problem which guided it towards high scoring priority programs.   But, second, such a description in fact appears \textsl{inapt} as the corpus $\mathcal{C}$ on which $M$ was trained had nothing to do with cap sets and the fact that the learning algorithm $\lambda$ which was employed was completely general.\footnote{\citet{Romera-Paredes2024} obtained their results using a particular large language model called Codey \citep{Google2023}.  This system was trained on human-generated computer code in high-level languages drawn from repositories such as GitHub.   This corpora included Python programs -- as well those in a number of languages -- which (one assumes) had been composed to carry out many disparate tasks.   But it appears that the training data had nothing  specifically to do with the cap set problem, combinatorics, or even mathematics in general.}   

In light of this, one is left to wonder exactly \textsl{why} the specific combination of techniques employed by \citep{Romera-Paredes2024} was successful.  The authors take some initial steps towards addressing this -- e.g. they report that they were unable to find the same admissible sets using a $\textsc{Sat}$-solver, genetic algorithms based on selection procedures involving random mutations, or a different large language model.\footnote{Nonetheless, the reported results are thus experimental in character and unsystematic in their failure to explore the full ranges of possible parameter values.  This illustrates how the questions highlighted above are connected to the general problem of \textsl{benchmarking} in artificial intelligence -- i.e. how can we infer from the behavior of a given system on a specific set of trials that its performance achieves a general standard or is superior to other systems in the general case?}  They additionally conjecture that FunSearch is successful because it has a bias for producing solutions  which ``in a loose sense $\ldots$ have low Kolmogorov complexity'' (p. 473) -- i.e. that the priority programs it generates tend to be small in size.\footnote{No attempt was made to precisify this claim nor is any evidence presented.  But even if it is true, it appears to offer light insight into why FunSearch was successful.  This is so for the simple reason that the maximally ``naive'' algorithm for finding an exact value for $\mathrm{cap}(8)$ by searching exhaustively $2^{3^8}$ possible subsets of $\{0,1,2\}^8$ to find the largest cap set itself admits to a very short program.  For this reason, it seems likely that an optimal algorithm for solving the problem efficiently would itself tend to be both complex (e.g. in virtue of being domain specific) rather than the simple Python priority programs by which FunSearch was used to obtain the (presumably) non-exact results in question. \label{kolmogorovnote}}  All of these claims await assessment.  But in addition to our programmatic observations about (C1) and (C2), they would presumably need to be addressed before we can hope to provide a satisfactory to the questions highlighted above.   

\section{Conclusion}

The argument of \S 3 concludes that proof discovery in mathematics is \textsl{inherently difficult} -- i.e. resistant to automation using techniques which can be carried out using computing machinery which can be constructed and and applied in practice.  On the other hand, the case studies reported in \S 4 have been repeatedly reported as breakthroughs for the application of artificial intelligence in mathematics.   But as we can now see, these are also of exactly the form which one would have predicted before the fact were amenable to automated resolution.   In particular, the statements obtained can all be formulated as statements of low logical complexity.   Thus not only do we known in advance that they are amenable to brute force search,  but this is in fact how the demonstrations in question proceeded.   

Brute force, is of course, also a method for which we have already ample evidence that computers outperform human mathematicians.  This is aptly illustrated by traditional ``computer proofs'' which do not involve techniques directly associated with artificial intelligence.   As such, our case studies also suggest that while the role of techniques like $\textsc{Sat}$-solvers and machine learning may come to play a larger role in mathematics, their influence will remain evolutionary rather than revolutionary.

We will now adduce two additional observations refining our prior claim that the examples of \S 4 are indeed \textsl{special cases} of proof discovery:\footnote{See \citep{Avigad2022, Avigad2024} for a further comparison of the results thus far obtained by automated methods and those which are regarded as significant within contemporary mathematical practice.}
\begin{itemize}
\item[(C$1'$)] Brute force search is not characteristic of how longstanding open questions in mathematics have historically been resolved.
\item[(C$2'$)] The logical complexity of many open questions is greater than $\Sigma^0_1$.   
\end{itemize}

An initial observation in regard to (C$2'$) is that the natural formalization of the  longstanding open questions consider in \S 2 which are related to number theory is at least $\Pi^0_1$ and often $\Pi^0_2$.   This is paradigmatically true of statements like the Goldbach and Twin Prime conjectures which can be directly expressed in language of first-order arithmetic.  But it also known to be true of statements like the Riemann Hypothesis or the $\mathbf{P} \neq \mathbf{NP}$ which can also be shown to be equivalent to statements of $\mathcal{L}_a$  as a consequence of additional lemmas and standard coding techniques.\footnote{\citet{Godel1931a} famously called attention to the fact that the Goldbach conjecture is naturally regarded as a $\Pi^0_1$-statement while \citet{Turing1939} made a similar observation about the formalization of the Riemann Hypothesis [RH] as a $\Pi^0_2$-statement.  This was then improved to $\Pi^0_1$ by \citet{Kreisel1952d,Kreisel1958} who explicitly constructed a primitive recursive predicate $B(x)$ such that $\forall x B(x)$ is provably equivalently to RH. (See also, e.g., \citep{Lagarias2002}.   The recent survey \citep{Broughan2023} describes efforts to reduce this to $\Delta^0_1$, which in turn would imply RH is either provable or refutable in Peano arithmetic.   But at the time of writing -- and to the best of our knowledge -- this hope remains unrealized.)  Statements which at present seem to be intrinsically $\Pi^0_2$ include $\mathbf{P} \neq \mathbf{NP}$ (see  \citealp{Ben-David1992}) and the Collatz Conjecture (see \citealp{Kurtz2007}).}

The fact that these questions have remained open -- sometimes for centuries -- suggest that we should at least not dismiss as naive the hypothesis that the logical complexity of mathematical statements provides a \textsl{prima facie} indication of their  discovermental complexity.    Of course this invites (at least) the following questions:

\vspace{.5ex}
\noindent \textsl{Is logical complexity -- e.g. as measured in the familiar manner of quantifier alternations -- an \textsl{intrinsic} feature of mathematical propositions?  If so, how do we determine the appropriate signature in which to formalize a proposition so that its complexity may be read off?  Are such assignments absolute or relative to a base theory over which potential complexity-reducing equivalences can be proven?\footnote{Suppose, for instance, that we discover that $\varphi$ is provable in $\mathsf{S}$.  Then (trivially) $\mathsf{S} \vdash_L \varphi \leftrightarrow \psi$ for $\psi$ any other statement provable in $\mathsf{S}$ -- e.g. a $\Delta^0_0$-formula like $2 + 2 = 4$.  A refinement of the prior question is thus as follows: To the extent which discovermental complexity can be regarded as intrinsic to a statement, is this complexity fixed by $\varphi$'s \textsl{surface form} (e.g. as reported by its initial statement in the literature), its \textsl{intermediate form} (e.g. once ``obvious'' transformations and lemmas have been taken into account), or its \textsl{deep$_{\mathsf{S}}$ form} (e.g. relative to ``non-obvious'' or even as yet undiscovered theorems derivable in $\mathsf{S}$)?} In the case of \textsl{non-arithmetical} -- e.g. geometrical, algebraic, analytic, or set-theoretic -- statements, does logical complexity still track discovermental  complexity? } 
\vspace{.5ex}

While these questions arise naturally, they appear not to have received answers.  But if we take for granted that at least some of the statements mentioned above are ``intrinsically'' $\Pi^0_1$ or higher, then a familiar observation is that they are least potentially undecidable in theories satisfying the hypotheses of Theorem 1.  A related point is that $\Pi^0_1$-statements are dual to $\Sigma^0_1$-statements in terms of the possibility of resolution by searching for a witness.   This is to say that while a $\Pi^0_1$-statement $\varphi \equiv \forall x \psi(x)$ can be \textsl{refuted} by conducting by finding a counterexample $n \in \mathbb{N}$ such that $\neg \psi(\overline{n})$, it cannot be \textsl{proven} in this manner.  This shows that the statements in question (if true) cannot be directly resolved by brute force.  

This in turn suggests that there may indeed be principled reasons to suspect that many open questions will remain outside the scope of the automated methods. But a qualification is required if we also take seriously the formalistic assumption which we employed in \S 2 to underpin premise (P2.ii) of the basic argument.  For if mathematical truth is  equated with derivability from a computable set of axioms, then the truth of \textsl{arbitrary} sentences $\varphi$ is reduced to the $\Sigma^0_1$-question of whether $\varphi \in \textsc{Prov}^S_{\Gamma}$.  And of course statements of this form \textsl{can} be resolved in the positive by the use of traditional automated theorem proving techniques to (essentially) perform a brute force search through the relevant domain of formal proofs.   

But significance of this observation appears limited in two respects.  First, one might deny the formalistic assumption while still accepting (P2.ii) in light of the detailed analysis of the relative complexity of  $\textsc{Prov}^S_{\Gamma}$ for different theories and logics as described in \S 3.  Second, the task of applying automated theorem proving techniques to statements more complex than those considered in \S 4 would seem to present considerable challenges.  While these cannot be considered in detail here, some of the difficulties can be appreciated by considering (C$1'$) directly.

An initial observation is that brute force is traditionally understood in contrast to other \textsl{non-brute} methods of computation and proof.   Such methods will characteristically \textsl{not} take the form of enumerative generate-and-test procedures, will \textsl{not} be confined to confirming $\Sigma^0_1$-statements,  and will prototypically be \textsl{domain specific}.  This falls short of providing an analysis of the relevant distinction. But even cursory examination suggests that non-brute methods are typically involved when longstanding open questions have historically been resolved -- e.g. by employing analogical reasoning, ``impure'' concepts and methods from \textsl{prima facie} unrelated domains, or by developing novel ``abstract'' theories or notions, etc.\footnote{Readers will likely be able to adduce their own favorite examples.  But some familiar illustrations include Galois's introduction of permutation groups in order to prove the non-existence of a quintic formula, the use of analytic or complex analytic methods in the original proofs of results like Dirichlet's Theorem or the Prime Number Theorem, the complicated architecture arising in the proof of Szemer\'edi's Theorem, and the introduction of the method of forcing in Cohen's proof of the independence of the Continuum Hypothesis.}

This does not rule out that automated methods might one day be applied to resolve a high-profile open problem in a brute force-like manner -- in fact it is exactly this which the third quote in the epigraph predicts.  But further consideration of how the methods in question currently operate in conjunction with logical observations about the complexity of formulas expressing the relevant statements suggests that if such a point ever arrives, it may still be a long way off.\footnote{For note that even the sort of sophisticated combination of automated theorem proving and machine learning described at the beginning of of \S 4.3 encounters problems discovering proofs which require instances of schema like mathematical induction which do not appear in the libraries on which they are trained (see, e.g., \citealp{Blanchette2016,Gauthier2023}).  But on the other hand, theories like $\mathsf{PA}$ and its extensions do not admit full cut elimination.  A simple consequence of this is that given a $\mathcal{L}_a$-sentence $\varphi$ of $\Pi^0_1$ or higher complexity, there is no effective means of determining an upper bound on the complexity of the instances of the induction scheme which might be required for a potential proof of $\varphi$.  Thus although provers of the sort in question are in principle capable of deriving statements of arbitrary complexity, they appear to be no better off than human mathematicians in terms of needing to find appropriate axioms to instantiate schema.   But at least as yet, such systems are unlike human mathematicians in the sense that they lack the facility to reason in the ``non-brute'' manner indicated above to find the relevant instances.} 

Of course the story of automated methods in mathematics will not end here.  In fact it seems all but certain that efforts to apply artificial intelligence-inspired technologies to proof discovery will continue to accelerate.  We have presented some reasons to think that such cases will continue to satisfy (C1) and (C2).   But we also hope that the foregoing has also illustrated the abundance of under-explored epistemological questions raised by both the limitative results surveyed in \S 3 and the positive case studies surveyed in \S 4.  In regard to the latter, we will conclude by returning to an intriguing aspect of the results described in \S 4.2 and \S 4.3.  

While we have argued that these are examples of brute force, we have seen that the relevant witnessing objects -- e.g. a certificate of unsatisfiability in the case of BPT or a priority program for generating cap sets -- were found after searching only a small fraction of the (astronomically large but still finite) space at issue. In fact, without this reduction it seems unlikely that the results could have been obtained with current or foreseeable technology.  But we are still drawn to ask:

\vspace{.5ex}

\noindent \textsl{Is it possible to explain the success in these cases in virtue of the formal properties of either the methods or the mathematical statements to which they were applied?   If so, does such an explanation provide any reason to believe that the methods will be successful when applied to other ``naturally occurring'' open problems in the relevant domains?    If not, is there any reason why we should not regard their success in the particular cases as mere ``coincidences'' or ``good luck''?}

\vspace{.5ex}

One way in these questions might be addressed is by providing a systematic account of the \textsl{division of labor} between human and automated systems -- e.g. in regard to finding traditional proofs which reduce the logical complexity of standing open problems or in fine-tuning of parameters of automated systems which are typically varied over the course of the \textsl{unsuccessful} application prior to a successful execution.    More generally, however, these questions further motivate the development of an account of what might be called \textsl{natural case complexity} as previewed in \S 2.3 -- e.g. one which weights formulas relative to the likelihood that we will formulate them explicitly or find them important within the conduct of our established mathematical practices.  \citet{Rabin1974} already called attention to the significance of this problem.  But we know of no general attempt to address it systematically.

{\footnotesize

\bibliographystyle{rsl}

\begin{thebibliography}{}

\bibitem[\protect\citeauthoryear{Aaronson}{Aaronson}{2016}]{Aaronson2016}
Aaronson, S. (2016).
\newblock $\mathbf{P} \stackrel{?}{=} \mathbf{NP}$.
\newblock In Nash, J. \& Rassias, M., editors, {\em Open problems in
  mathematics}, pp.\  1--122. Berlin: Springer.

\bibitem[\protect\citeauthoryear{Adem}{Adem}{2024}]{Adem2024}
Adem, A., editor (2024).
\newblock {\em Bulletin of the American Mathematical Society}, Volume 61 (2).

\bibitem[\protect\citeauthoryear{Appel, \& Haken}{Appel \&
  Haken}{1977}]{Appel1977a}
Appel, K., \& Haken, W. (1977).
\newblock {Every planar map is four colorable. Part I: Discharging}.
\newblock {\em Illinois Journal of Mathematics\/}~{\bf 21\/}(3), 429--490.

\bibitem[\protect\citeauthoryear{Arana}{Arana}{2015}]{Arana2015}
Arana, A. (2015).
\newblock {On the Depth of Szemer{\'e}di's Theorem}.
\newblock {\em Philosophia Mathematica\/}~{\bf 23\/}(2), 163--176.

\bibitem[\protect\citeauthoryear{Arana, \& Stafford}{Arana \&
  Stafford}{2023}]{Arana2023}
Arana, A., \& Stafford, W. (2023).
\newblock On the difficulty of discovering mathematical proofs.
\newblock {\em Synthese\/}~{\bf 202\/}(2), 38.

\bibitem[\protect\citeauthoryear{Avigad}{Avigad}{2021}]{Avigad2021}
Avigad, J. (2021).
\newblock Reliability of mathematical inference.
\newblock {\em Synthese\/}~{\bf 198\/}(8), 7377--7399.

\bibitem[\protect\citeauthoryear{Avigad}{Avigad}{2022}]{Avigad2022}
Avigad, J. (2022).
\newblock Varieties of mathematical understanding.
\newblock {\em Bulletin of the American Mathematical Society\/}~{\bf 59\/}(1),
  99--117.

\bibitem[\protect\citeauthoryear{Avigad}{Avigad}{2024}]{Avigad2024}
Avigad, J. (2024).
\newblock Automated reasoning for mathematics.
\newblock In Benzm{\"u}ller, C., Heule, M.~J., \& Schmidt, R.~A., editors, {\em
  Automated Reasoning}, pp.\  3--20.

\bibitem[\protect\citeauthoryear{Ben-David, \& Halevi}{Ben-David \&
  Halevi}{1992}]{Ben-David1992}
Ben-David, S., \& Halevi, S. (1992).
\newblock On the independence of {P} versus {NP}.
\newblock Technical Report 714, Technion.

\bibitem[\protect\citeauthoryear{Berman}{Berman}{1980}]{Berman1980}
Berman, L. (1980).
\newblock The complexity of logical theories.
\newblock {\em Theoretical Computer Science\/}~{\bf 11\/}(1), 71--77.

\bibitem[\protect\citeauthoryear{Biere, Heule, van Maaren, \& Walsh}{Biere
  et~al.}{2021}]{Biere2021}
Biere, A., Heule, M., van Maaren, H., \& Walsh, T. (2021).
\newblock {\em Handbook of satisfiability}, Volume 336.
\newblock IOS press.

\bibitem[\protect\citeauthoryear{Blanchette, Kaliszyk, Paulson, \&
  Urban}{Blanchette et~al.}{2016}]{Blanchette2016}
Blanchette, J.~C., Kaliszyk, C., Paulson, L.~C., \& Urban, J. (2016).
\newblock {Hammering towards QED}.
\newblock {\em Journal of Formalized Reasoning\/}~{\bf 9\/}(1), 101--148.

\bibitem[\protect\citeauthoryear{Borwein, \& Bailey}{Borwein \&
  Bailey}{2008}]{Borwein2008}
Borwein, J., \& Bailey, D. (2008).
\newblock {\em Mathematics by experiment}.
\newblock CRC Press Boca Raton, FL, USA.

\bibitem[\protect\citeauthoryear{Brakensiek, Heule, Mackey, \&
  Narv{\'a}ez}{Brakensiek et~al.}{2020}]{Brakensiek2020}
Brakensiek, J., Heule, M., Mackey, J., \& Narv{\'a}ez, D. (2020).
\newblock {The resolution of Keller's conjecture}.
\newblock In {\em International Joint Conference on Automated Reasoning}, pp.\
  48--65. Springer.

\bibitem[\protect\citeauthoryear{Broughan}{Broughan}{2023}]{Broughan2023}
Broughan, K. (2023).
\newblock {\em Equivalents of the Riemann Hypothesis: Volume 3, Further Steps
  towards Resolving the Riemann Hypothesis}.
\newblock Cambridge University Press.

\bibitem[\protect\citeauthoryear{Buss}{Buss}{1998}]{Buss1998b}
Buss, S.~R. (1998).
\newblock First-{O}rder {P}roof {T}heory of {A}rithmetic.
\newblock In {\em Handbook of {P}roof {T}heory}, Volume 137, pp.\  79--147.
  Amsterdam: North-Holland.

\bibitem[\protect\citeauthoryear{Carlson, Jaffe, \& Wiles}{Carlson
  et~al.}{2006}]{Carlson2006}
Carlson, J., Jaffe, A., \& Wiles, A., editors (2006).
\newblock {\em {The Millennium Prize Problems}}.
\newblock Providence, R. I.: American Mathematical Society.

\bibitem[\protect\citeauthoryear{Castelvecchi}{Castelvecchi}{2024}]{Castelvecchi2024}
Castelvecchi, D. (2024).
\newblock {DeepMind AI outdoes human mathematicians on unsolved problem}.
\newblock {\em Nature\/}~{\bf 625\/}(7993), 12--13.

\bibitem[\protect\citeauthoryear{Cook, \& Reckhow}{Cook \&
  Reckhow}{1979}]{Cook1979}
Cook, S., \& Reckhow, R. (1979).
\newblock The relative efficiency of propositional proof systems.
\newblock {\em The Journal of Symbolic Logic\/}~{\bf 44\/}(01), 36--50.

\bibitem[\protect\citeauthoryear{Cooper}{Cooper}{2004}]{Cooper2004}
Cooper, S. (2004).
\newblock {\em Computability Theory}.
\newblock Boca Raton: Chapman \& Hall.

\bibitem[\protect\citeauthoryear{Cormen, Leiserson, \& Rivest}{Cormen
  et~al.}{2005}]{Cormen2005}
Cormen, T., Leiserson, C., \& Rivest, R. (2005).
\newblock {\em Introduction to algorithms\/} (Second ed.).
\newblock MIT Press.

\bibitem[\protect\citeauthoryear{Crandall, \& Pomerance}{Crandall \&
  Pomerance}{2005}]{Crandall2005}
Crandall, R., \& Pomerance, C. (2005).
\newblock {\em Prime {N}umbers\/} (Second ed.).
\newblock New York: Springer.

\bibitem[\protect\citeauthoryear{Davies, Veli{\v{c}}kovi{\'c}, Buesing,
  Blackwell, Zheng, Toma{\v{s}}ev, Tanburn, Battaglia, Blundell, Juh{\'a}sz,
  et~al.}{Davies et~al.}{2021}]{Davies2021}
Davies, A., Veli{\v{c}}kovi{\'c}, P., Buesing, L., Blackwell, S., Zheng, D.,
  Toma{\v{s}}ev, N., Tanburn, R., Battaglia, P., Blundell, C., Juh{\'a}sz, A.,
  et~al. (2021).
\newblock {Advancing mathematics by guiding human intuition with AI}.
\newblock {\em Nature\/}~{\bf 600\/}(7887), 70--74.

\bibitem[\protect\citeauthoryear{Davis}{Davis}{1957}]{Davis1957}
Davis, M. (1957).
\newblock {A program for Presburger's algorithm}.
\newblock In {\em Summer Institute for Symbolic Logic}, pp.\  215--227. Cornell
  Mathematics Department.

\bibitem[\protect\citeauthoryear{Dean}{Dean}{2019}]{Dean2019c}
Dean, W. (2019).
\newblock Computational complexity theory and the philosophy of mathematics.
\newblock {\em Philosophia Mathematica\/}~{\bf 27\/}(3), 381--439.

\bibitem[\protect\citeauthoryear{Dershowitz, \& Jouannaud}{Dershowitz \&
  Jouannaud}{1990}]{Dershowitz1990}
Dershowitz, N., \& Jouannaud, J.-P. (1990).
\newblock Rewrite systems.
\newblock In Van~Leeuwen, J., editor, {\em {Handbook of theoretical computer
  science}}, Volume~B, pp.\  243--320. Elsevier.

\bibitem[\protect\citeauthoryear{Detlefsen}{Detlefsen}{1990}]{Detlefsen1990}
Detlefsen, M. (1990).
\newblock {On an alleged refutation of Hilbert's program using G{\"o}del's
  first incompleteness theorem}.
\newblock {\em Journal of Philosophical Logic\/}~{\bf 19\/}(4), 343--377.

\bibitem[\protect\citeauthoryear{Detlefsen}{Detlefsen}{1996}]{Detlefsen1996}
Detlefsen, M. (1996).
\newblock Philosophy of {M}athematics in the {T}wentieth {C}entury.
\newblock In Shanker, S.~G., editor, {\em Philosophy of {S}cience, {L}ogic and
  {M}athematics in the {T}wentieth {C}entury}, Volume~9 of {\em Routledge
  History of Philosophy}, pp.\  50--123. London: Routledge.

\bibitem[\protect\citeauthoryear{Detlefsen, \& Luker}{Detlefsen \&
  Luker}{1980}]{Detlefsen1980}
Detlefsen, M., \& Luker, M. (1980).
\newblock The four-color theorem and mathematical proof.
\newblock {\em The Journal of Philosophy\/}~{\bf 77}, 803--802.

\bibitem[\protect\citeauthoryear{Edel}{Edel}{2004}]{Edel2004}
Edel, Y. (2004).
\newblock Extensions of generalized product caps.
\newblock {\em Designs, Codes and Cryptography\/}~{\bf 31}, 5--14.

\bibitem[\protect\citeauthoryear{Ellenberg, \& Gijswijt}{Ellenberg \&
  Gijswijt}{2017}]{Ellenberg2017}
Ellenberg, J.~S., \& Gijswijt, D. (2017).
\newblock On large subsets of with no three-term arithmetic progression.
\newblock {\em Annals of Mathematics\/}~{\bf 186}, 339--343.

\bibitem[\protect\citeauthoryear{Ewald}{Ewald}{1996}]{Ewald1996}
Ewald, W. (1996).
\newblock {\em From {K}ant to {H}ilbert: {A} {S}ource {B}ook in the
  {F}oundations of {M}athematics.}
\newblock New York: Oxford University Press.

\bibitem[\protect\citeauthoryear{Ewald, \& Sieg}{Ewald \&
  Sieg}{2013}]{Hilbert2013}
Ewald, W., \& Sieg, W., editors (2013).
\newblock {\em {David Hilbert's Lectures on the Foundations of Logic and
  Arithmetic 1917 -- 1933}}.
\newblock Berlin: Springer.

\bibitem[\protect\citeauthoryear{Feferman, et~al.}{Feferman
  et~al.}{1986}]{Godel1986}
Feferman, S., et~al., editors (1986).
\newblock {\em {Kurt G\"odel Collected Works. Vol. I. Publications 1929-1936}}.
\newblock Oxford: Oxford Univeristy Press.

\bibitem[\protect\citeauthoryear{Fischer, \& Rabin}{Fischer \&
  Rabin}{1974}]{Fischer1974}
Fischer, M., \& Rabin, M. (1974).
\newblock {Super-exponential complexity of Presburger arithmetic}.
\newblock In Karp, R., editor, {\em Complexity of computation}, {SIAM-AMS
  Symposium in Applied Mathematics}.

\bibitem[\protect\citeauthoryear{Gandy}{Gandy}{1980}]{Gandy1980}
Gandy, R. (1980).
\newblock {Church's Thesis and principles for mechanisms}.
\newblock In J.~Barwise, H.~K. \& Kunen, K., editors, {\em The Kleene
  Symposium}, Volume 101, pp.\  123--148. North Holland.

\bibitem[\protect\citeauthoryear{Gauthier, Brown, Janota, \& Urban}{Gauthier
  et~al.}{2023}]{Gauthier2023}
Gauthier, T., Brown, C.~E., Janota, M., \& Urban, J. (2023).
\newblock A mathematical benchmark for inductive theorem provers.
\newblock arXiv preprint arXiv:2304.02986.

\bibitem[\protect\citeauthoryear{Gauthier, Kaliszyk, \& Urban}{Gauthier
  et~al.}{2016}]{Gauthier2016}
Gauthier, T., Kaliszyk, C., \& Urban, J. (2016).
\newblock Initial experiments with statistical conjecturing over large formal
  corpora.
\newblock In et~al., A.~K., editor, {\em Joint Proceedings of the FM4M, MathUI,
  and ThEdu Workshops, Doctoral Program, and Work in Progress at the Conference
  on Intelligent Computer Mathematics}, Volume 1785, pp.\  219--228.

\bibitem[\protect\citeauthoryear{Gauthier, Ol{\v{s}}{\'a}k, \& Urban}{Gauthier
  et~al.}{2023}]{Gauthier2023a}
Gauthier, T., Ol{\v{s}}{\'a}k, M., \& Urban, J. (2023).
\newblock Alien coding.
\newblock {\em International Journal of Approximate Reasoning\/}~{\bf
  162\/}(109009), 1--26.

\bibitem[\protect\citeauthoryear{Gauthier, \& Urban}{Gauthier \&
  Urban}{2023}]{Gauthier2023b}
Gauthier, T., \& Urban, J. (2023).
\newblock Learning program synthesis for integer sequences from scratch.
\newblock {\em Proceedings of the AAAI Conference on Artificial
  Intelligence\/}~{\bf 37\/}(6), 7670--7677.

\bibitem[\protect\citeauthoryear{G\"odel}{G\"odel}{1931}]{Godel1931a}
G\"odel, K. (1931).
\newblock {On} formally undecidable propositions of {{\em {Principia
  Mathematica}\/}} and related systems~{I}.
\newblock Reprinted in \cite{Godel1986}, pp. 144-195.

\bibitem[\protect\citeauthoryear{Gonthier, Asperti, Avigad, Bertot, Cohen,
  Garillot, Le~Roux, Mahboubi, O'Connor, Biha, et~al.}{Gonthier
  et~al.}{2013}]{Gonthier2013}
Gonthier, G., Asperti, A., Avigad, J., Bertot, Y., Cohen, C., Garillot, F.,
  Le~Roux, S., Mahboubi, A., O'Connor, R., Biha, S.~O., et~al. (2013).
\newblock A machine-checked proof of the odd order theorem.
\newblock In {\em International conference on interactive theorem proving},
  pp.\  163--179. Springer.

\bibitem[\protect\citeauthoryear{Google}{Google}{2023}]{Google2023}
Google (2023).
\newblock Code models overview.
\newblock
  \url{https://cloud.google.com/vertex-ai/generative-ai/docs/code/code-models-overview}.

\bibitem[\protect\citeauthoryear{Greenlaw, Hoover, \& Ruzzo}{Greenlaw
  et~al.}{1995}]{Greenlaw1995}
Greenlaw, R., Hoover, H., \& Ruzzo, W. (1995).
\newblock {\em {Limits to parallel computation: \textbf{P}-completeness
  theory}}.
\newblock Oxford University Press, USA.

\bibitem[\protect\citeauthoryear{Grochow}{Grochow}{2019}]{Grochow2019}
Grochow, J. (2019).
\newblock {New applications of the polynomial method: The cap set conjecture
  and beyond}.
\newblock {\em Bulletin of the American Mathematical Society\/}~{\bf 56\/}(1),
  29--64.

\bibitem[\protect\citeauthoryear{Haase}{Haase}{2018}]{Haase2018}
Haase, C. (2018).
\newblock {A survival guide to Presburger arithmetic}.
\newblock {\em ACM SIGLOG News\/}~{\bf 5\/}(3), 67--82.

\bibitem[\protect\citeauthoryear{Hales, Adams, Bauer, Dang, Harrison,
  Le~Truong, Kaliszyk, Magron, McLaughlin, Nguyen, et~al.}{Hales
  et~al.}{2017}]{Hales2017}
Hales, T., Adams, M., Bauer, G., Dang, T.~D., Harrison, J., Le~Truong, H.,
  Kaliszyk, C., Magron, V., McLaughlin, S., Nguyen, T.~T., et~al. (2017).
\newblock {A formal proof of the Kepler conjecture}.
\newblock In {\em Forum of mathematics, Pi}, Volume~5. Cambridge University
  Press.

\bibitem[\protect\citeauthoryear{Hales}{Hales}{2014}]{Hales2014}
Hales, T.~C. (2014).
\newblock Developments in formal proofs.
\newblock {\em S{\'e}minaire Bourbaki\/}~{\bf 66\/}(1086), 1--23.

\bibitem[\protect\citeauthoryear{Hardy}{Hardy}{1940}]{Hardy1940}
Hardy, G. (1940).
\newblock {\em A mathematician's apology}.
\newblock Cambridge University Press.

\bibitem[\protect\citeauthoryear{Hardy}{Hardy}{1929}]{Hardy1929}
Hardy, G.~H. (1929).
\newblock Mathematical proof.
\newblock {\em Mind\/}~{\bf 38\/}(149), 1--25.

\bibitem[\protect\citeauthoryear{Henkin, Tarski, \& Monk}{Henkin
  et~al.}{1971}]{Henkin1971}
Henkin, L., Tarski, A., \& Monk, D. (1971).
\newblock {\em Cylindric algebras}.
\newblock North Holland.

\bibitem[\protect\citeauthoryear{Hennessy, Patterson, \& Goldberg}{Hennessy
  et~al.}{2003}]{Hennessy2003}
Hennessy, J., Patterson, D., \& Goldberg, D. (2003).
\newblock {\em {Computer architecture: a quantitative approach}}.
\newblock Morgan Kaufmann.

\bibitem[\protect\citeauthoryear{Heule, Kullmann, \& Marek}{Heule
  et~al.}{2016}]{Heule2016}
Heule, M., Kullmann, O., \& Marek, V.~W. (2016).
\newblock {Solving and Verifying the Boolean Pythagorean Triples Problem via
  Cube-and-Conquer}.
\newblock In {\em International Conference on Theory and Applications of
  Satisfiability Testing}, pp.\  228--245. Springer.

\bibitem[\protect\citeauthoryear{Heule, \& Kullmann}{Heule \&
  Kullmann}{2017}]{Heule2017}
Heule, M.~J., \& Kullmann, O. (2017).
\newblock The science of brute force.
\newblock {\em Communications of the ACM\/}~{\bf 60\/}(8), 70--79.

\bibitem[\protect\citeauthoryear{Heule, \& Scheucher}{Heule \&
  Scheucher}{2024}]{Heule2024}
Heule, M.~J., \& Scheucher, M. (2024).
\newblock Happy ending: An empty hexagon in every set of 30 points.
\newblock In {\em International Conference on Tools and Algorithms for the
  Construction and Analysis of Systems}, pp.\  61--80. Springer.

\bibitem[\protect\citeauthoryear{Hilbert}{Hilbert}{1900}]{Hilbert1900}
Hilbert, D. (1900).
\newblock Mathematische {P}robleme. {V}ortrag, gehalten auf dem internationalen
  {M}athematiker-{C}ongress zu {P}aris 1900.
\newblock In {\em Nachrichten von der Gesellschaft der Wissenschaften zu
  G\"ottingen, Mathematisch-Physikalische Klasse}, pp.\  253--297.
\newblock English transalation as ``Mathematical Problems'' in
  \cite{Ewald1996}, pp. 1096-1105.

\bibitem[\protect\citeauthoryear{Hilbert}{Hilbert}{1930}]{Hilbert1930a}
Hilbert, D. (1930).
\newblock Naturerkennen und logik.
\newblock {\em Die Naturwissenschaften\/}~{\bf 18}, 959--63.
\newblock Translated by W. Ewald as ``Logic and the knowledge of nature'' and
  republshed in \citep{Ewald1996}, Volume II.

\bibitem[\protect\citeauthoryear{Hilbert, \& Ackermann}{Hilbert \&
  Ackermann}{1928}]{Hilbert1928}
Hilbert, D., \& Ackermann, W. (1928).
\newblock {\em Grundz{\"u}ge der theoretischen Logik\/} (First ed.).
\newblock Berlin: Springer.
\newblock Reprinted in \cite{Hilbert2013}.

\bibitem[\protect\citeauthoryear{Holland}{Holland}{1992}]{Holland1992}
Holland, J. (1992).
\newblock {\em Adaptation in natural and artificial systems: an introductory
  analysis with applications to biology, control, and artificial intelligence}.
\newblock MIT press.

\bibitem[\protect\citeauthoryear{Huntington}{Huntington}{1933}]{Huntington1933}
Huntington, E.~V. (1933).
\newblock {New sets of independent postulates for the algebra of logic, with
  special reference to Whitehead and Russell's \textsl{Principia mathematica}}.
\newblock {\em Transactions of the American Mathematical Society\/}~{\bf
  35\/}(1), 274--304.

\bibitem[\protect\citeauthoryear{Kinyon, Veroff, \& Petr}{Kinyon
  et~al.}{2013}]{Kinyon2013}
Kinyon, M., Veroff, R., \& Petr, V. (2013).
\newblock Loops with abelian inner mapping groups: An application of automated
  deduction.
\newblock In {\em Automated Reasoning and Mathematics}, pp.\  151--164.
  Springer.

\bibitem[\protect\citeauthoryear{Kolata}{Kolata}{1996}]{Kolata1996}
Kolata, G. (1996).
\newblock With major math proof, brute computers show flash of reasoning power.
\newblock {\em New York Times\/}.
\newblock Decemeber 10.

\bibitem[\protect\citeauthoryear{Kreisel}{Kreisel}{1952}]{Kreisel1952d}
Kreisel, G. (1952).
\newblock Some elementary inequalities.
\newblock {\em Indagationes Mathematicae\/}~{\bf 14}, 334--338.

\bibitem[\protect\citeauthoryear{Kreisel}{Kreisel}{1958}]{Kreisel1958}
Kreisel, G. (1958).
\newblock Mathematical significance of consistency proofs.
\newblock {\em The Journal of Symbolic Logic\/}~{\bf 23\/}(2), 155--182.

\bibitem[\protect\citeauthoryear{Kurtz, \& Simon}{Kurtz \&
  Simon}{2007}]{Kurtz2007}
Kurtz, S.~A., \& Simon, J. (2007).
\newblock {The undecidability of the generalized Collatz problem}.
\newblock In {\em International Conference on Theory and Applications of Models
  of Computation}, pp.\  542--553. Springer.

\bibitem[\protect\citeauthoryear{Lagarias}{Lagarias}{2002}]{Lagarias2002}
Lagarias, J.~C. (2002).
\newblock {An elementary problem equivalent to the Riemann hypothesis}.
\newblock {\em The American mathematical monthly\/}~{\bf 109\/}(6), 534--543.

\bibitem[\protect\citeauthoryear{Landau}{Landau}{1912}]{Landau1912}
Landau, E. (1912).
\newblock Gel{\"o}ste und ungel{\"o}ste probleme aus der theorie der
  primzahlverteilung und der riemannschen zetafunktion.
\newblock {\em Jahresbericht der Deutschen Mathematiker-Vereinigung\/}~{\bf
  21}, 208--228.

\bibitem[\protect\citeauthoryear{Lighthill}{Lighthill}{1973}]{Lighthill1973}
Lighthill, J. (1973).
\newblock Artificial intelligence: a general survey.
\newblock In {\em Artificial Intelligence: a paper symposium}. Science Research
  Council.

\bibitem[\protect\citeauthoryear{Lin, Cao, Huang, Yang, Liu, Li, \& Liang}{Lin
  et~al.}{2024}]{Lin2024}
Lin, X., Cao, Q., Huang, Y., Yang, Z., Liu, Z., Li, Z., \& Liang, X. (2024).
\newblock Atg: Benchmarking automated theorem generation for generative
  language models.
\newblock In {\em Findings of the Association for Computational Linguistics:
  NAACL 2024}, pp.\  4465--4480.

\bibitem[\protect\citeauthoryear{Lucas}{Lucas}{1961}]{Lucas1961}
Lucas, J.~R. (1961).
\newblock {Minds, Machines and G{\"o}del}.
\newblock {\em Philosophy\/}~{\bf 36\/}(137), 112--127.

\bibitem[\protect\citeauthoryear{MacKenzie}{MacKenzie}{2004}]{MacKenzie2001}
MacKenzie, D.~A. (2004).
\newblock {\em Mechanizing proof: computing, risk, and trust}.
\newblock MIT Press.

\bibitem[\protect\citeauthoryear{McCune}{McCune}{1997}]{McCune1997}
McCune, W. (1997).
\newblock Solution of the robbins problem.
\newblock {\em Journal of Automated Reasoning\/}~{\bf 19\/}(3), 263--276.

\bibitem[\protect\citeauthoryear{McMahon, Gordon, Gordon, \& Gordon}{McMahon
  et~al.}{2019}]{McMahon2019}
McMahon, L., Gordon, G., Gordon, H., \& Gordon, R. (2019).
\newblock {\em The joy of SET: The many mathematical dimensions of a seemingly
  simple card game}.
\newblock Princeton University Press.

\bibitem[\protect\citeauthoryear{Nash, \& Rassias}{Nash \&
  Rassias}{2016}]{Nash2016}
Nash, J.~F., \& Rassias, M.~T. (2016).
\newblock {\em Open problems in mathematics}.
\newblock Springer.

\bibitem[\protect\citeauthoryear{Newell, Shaw, \& Simon}{Newell
  et~al.}{1957}]{Newell1957}
Newell, A., Shaw, J.~C., \& Simon, H.~A. (1957).
\newblock Empirical explorations of the logic theory machine: a case study in
  heuristic.
\newblock In {\em Papers presented at the February 26-28, 1957, western joint
  computer conference: Techniques for reliability}, pp.\  218--230.

\bibitem[\protect\citeauthoryear{Parikh}{Parikh}{1971}]{Parikh1971}
Parikh, R. (1971).
\newblock Existence and feasibility in arithmetic.
\newblock {\em Journal of Symbolic Logic\/}~{\bf 36\/}(3), 494--508.

\bibitem[\protect\citeauthoryear{Polya}{Polya}{1954}]{Polya1954}
Polya, G. (1954).
\newblock {\em Mathematics and plausible reasoning: Induction and analogy in
  mathematics}, Volume~1.
\newblock Princeton University Press.

\bibitem[\protect\citeauthoryear{Post}{Post}{1944}]{Post1944}
Post, E. (1944).
\newblock Recursively enummerable sets of postive integers and their decision
  problems.
\newblock {\em Bulletin of the American Mathematical Society\/}~{\bf 50},
  284--316.

\bibitem[\protect\citeauthoryear{Presburger}{Presburger}{1930}]{Presburger1930}
Presburger, M. (1930).
\newblock \"{U}ber die {V}ollst\"andigkeit eines gewissen {S}ystems der
  {A}rithmetik ganzer {Z}ahlen, in welchem die {A}ddition als einzige
  {O}peration hervortritt.
\newblock In {\em Sprawozdanie z 1. Kongresu matematyk{\'o}w kraj{\'o}w
  slowianskich. Comptes-rendus du 1. Congr{\`e}s des math{\'e}maticiens des
  pays slaves, Warszawa, 1929}, pp.\  92--101. Warszawa: Ksiaznica.

\bibitem[\protect\citeauthoryear{Quaife}{Quaife}{1992}]{Quaife1992}
Quaife, A. (1992).
\newblock {Automated deduction in von Neumann-Bernays-G{\"o}del set theory}.
\newblock {\em Journal of Automated Reasoning\/}~{\bf 8\/}(1), 91--147.

\bibitem[\protect\citeauthoryear{Rabin}{Rabin}{1974}]{Rabin1974}
Rabin, M.~O. (1974).
\newblock Theoretical impediments to artificial intelligence.
\newblock In {\em Proceedings of the 6th IFIP Congress 1974, Stockholm}, pp.\
  615--619.

\bibitem[\protect\citeauthoryear{Roberts}{Roberts}{2023}]{Roberts2023}
Roberts, S. (2023).
\newblock {A.I. Is Coming for Mathematics, Too}.
\newblock {\em The New York Times Iternational Edition\/}.
\newblock July 18.

\bibitem[\protect\citeauthoryear{Robinson}{Robinson}{1965}]{Robinson1965d}
Robinson, J.~A. (1965).
\newblock A machine-oriented logic based on the resolution principle.
\newblock {\em Journal of the ACM\/}~{\bf 12\/}(1), 23--41.

\bibitem[\protect\citeauthoryear{Romera-Paredes, Barekatain, Novikov, Balog,
  Kumar, Dupont, Ruiz, Ellenberg, Wang, Fawzi, Kohli, \& Fawzi}{Romera-Paredes
  et~al.}{2024}]{Romera-Paredes2024}
Romera-Paredes, B., Barekatain, M., Novikov, A., Balog, M., Kumar, M.~P.,
  Dupont, E., Ruiz, F. J.~R., Ellenberg, J.~S., Wang, P., Fawzi, O., Kohli, P.,
  \& Fawzi, A. (2024).
\newblock Mathematical discoveries from program search with large language
  models.
\newblock {\em Nature\/}~{\bf 625\/}(7995), 468--475.

\bibitem[\protect\citeauthoryear{Rosser}{Rosser}{1936}]{Rosser1936}
Rosser, B. (1936).
\newblock {Extensions of some theorems of G{\"o}del and Church}.
\newblock {\em Journal of symbolic logic\/}~{\bf 1\/}(3), 87--91.

\bibitem[\protect\citeauthoryear{Russell, \& Norvig}{Russell \&
  Norvig}{2021}]{Russell2021}
Russell, S.~J., \& Norvig, P. (2021).
\newblock {\em Artificial intelligence: a modern approach\/} (4th ed.).
\newblock Pearson.

\bibitem[\protect\citeauthoryear{Sample}{Sample}{2023}]{Sample2023}
Sample, I. (2023).
\newblock {AI scientists make `exciting' discovery using chatbots to solve
  maths problems}.
\newblock {\em {The Guardian}\/}.
\newblock December 14.

\bibitem[\protect\citeauthoryear{Savitch, \& Stimson}{Savitch \&
  Stimson}{1979}]{Savitch1979}
Savitch, W., \& Stimson, M. (1979).
\newblock {Time bounded random access machines with parallel processing}.
\newblock {\em Journal of the ACM\/}~{\bf 26\/}(1), 103--118.

\bibitem[\protect\citeauthoryear{Shoenfield}{Shoenfield}{1967}]{Shoenfield1967}
Shoenfield, J. (1967).
\newblock {\em Mathematical logic}, Volume~21.
\newblock Addison-Wesley Reading.

\bibitem[\protect\citeauthoryear{Sieg}{Sieg}{2009}]{Sieg2009}
Sieg, W. (2009).
\newblock On computability.
\newblock In Irvine, A., editor, {\em Philosophy of mathematics}, Volume~4 of
  {\em Handbook of the philosophy of science}, pp.\  535--630. Amsterdam: North
  Holland.

\bibitem[\protect\citeauthoryear{Simon, \& Newell}{Simon \&
  Newell}{1958}]{Simon1958}
Simon, H.~A., \& Newell, A. (1958).
\newblock Heuristic problem solving: The next advance in operations research.
\newblock {\em Operations research\/}~{\bf 6\/}(1), 1--10.

\bibitem[\protect\citeauthoryear{Slaman}{Slaman}{1998}]{Slaman1998}
Slaman, T.~A. (1998).
\newblock Mathematical {D}efinability.
\newblock In Dales, H. \& Oliveri, G., editors, {\em Truth in {M}athematics},
  pp.\  233--251. New York: Oxford University Press.

\bibitem[\protect\citeauthoryear{Sloane}{Sloane}{2007}]{Sloane2007}
Sloane, N. (2007).
\newblock The on-line encyclopedia of integer sequences.
\newblock In Kauers, M., Kerber, M., Miner, R., \& Windsteiger, W., editors,
  {\em Towards Mechanized Mathematical Assistants: 14th Symposium, Calculemus
  2007, 6th International Conference, MKM 2007, Hagenberg, Austria, June 27-30,
  2007. Proceedings}, pp.\  130--130. Springer.

\bibitem[\protect\citeauthoryear{Smale}{Smale}{1998}]{Smale1998}
Smale, S. (1998).
\newblock Mathematical problems for the next century.
\newblock {\em The mathematical intelligencer\/}~{\bf 20\/}(2), 7--15.

\bibitem[\protect\citeauthoryear{Smorynski}{Smorynski}{1991}]{Smorynski1991a}
Smorynski, C. (1991).
\newblock {\em Logical number theory I: An introduction}.
\newblock Springer.

\bibitem[\protect\citeauthoryear{Smory{\'n}ski}{Smory{\'n}ski}{2020}]{Smorynski2020}
Smory{\'n}ski, C. (2020).
\newblock {\em Mathematical problems: An essay on their nature and importance}.
\newblock Springer.

\bibitem[\protect\citeauthoryear{Subercaseaux, \& Heule}{Subercaseaux \&
  Heule}{2023}]{Subercaseaux2023}
Subercaseaux, B., \& Heule, M. J.~H. (2023).
\newblock The packing chromatic number of the infinite square grid is 15.
\newblock In Sankaranarayanan, S. \& Sharygina, N., editors, {\em Tools and
  Algorithms for the Construction and Analysis of Systems}, pp.\  389--406.

\bibitem[\protect\citeauthoryear{Tarski}{Tarski}{1959}]{Tarski1959a}
Tarski, A. (1959).
\newblock What is {E}lementary {G}eometry?
\newblock In {\em The Axiomatic Method.}, Studies in Logic and the Foundations
  of Mathematics, pp.\  16--29. Amsterdam: North-Holland.

\bibitem[\protect\citeauthoryear{Tarski, Mostowski, \& Robinson}{Tarski
  et~al.}{1953}]{Tarski1953}
Tarski, A., Mostowski, A., \& Robinson, R. (1953).
\newblock {\em Undecidable {T}heories}.
\newblock Studies in Logic and the Foundations of Mathematics. Amsterdam:
  North-Holland.

\bibitem[\protect\citeauthoryear{Trakhtenbrot}{Trakhtenbrot}{1984}]{Trakhtenbrot1984}
Trakhtenbrot, B. (1984).
\newblock {A survey of Russian approaches to \textsl{perebor} (brute-force
  searches) algorithms}.
\newblock {\em Annals of the History of Computing\/}~{\bf 6\/}(4), 384--400.

\bibitem[\protect\citeauthoryear{Troelstra, \& Schwichtenberg}{Troelstra \&
  Schwichtenberg}{2000}]{Troelstra2000}
Troelstra, A., \& Schwichtenberg, H. (2000).
\newblock {\em Basic Proof Theory\/} (Second ed.).
\newblock Cambridge University Press.

\bibitem[\protect\citeauthoryear{Turing}{Turing}{1939}]{Turing1939}
Turing, A. (1939).
\newblock Systems of logic based on ordinals.
\newblock {\em Proceedings of the London Mathematical Society\/}~{\bf 2\/}(1),
  161--228.

\bibitem[\protect\citeauthoryear{Tymoczko}{Tymoczko}{1979}]{Tymoczko1979}
Tymoczko, T. (1979).
\newblock The four-color problem and its philosophical significance.
\newblock {\em Journal of Philosophy\/}~{\bf 76\/}(2), 57--83.

\bibitem[\protect\citeauthoryear{Tyrrell}{Tyrrell}{2023}]{Tyrrell2023}
Tyrrell, F. (2023).
\newblock New lower bounds for cap sets.
\newblock {\em Discrete analysis\/}~{\bf 20}, 1--18.

\bibitem[\protect\citeauthoryear{Urban, \& Vysko{\v{c}}il}{Urban \&
  Vysko{\v{c}}il}{2013}]{Urban2013}
Urban, J., \& Vysko{\v{c}}il, J. (2013).
\newblock {Theorem proving in large formal mathematics as an emerging AI
  field}.
\newblock In {\em Automated Reasoning and Mathematics}, pp.\  240--257.
  Springer.

\bibitem[\protect\citeauthoryear{van Emde~Boas}{van Emde~Boas}{1990}]{Boas1990}
van Emde~Boas, P. (1990).
\newblock {Machine models and simulations}.
\newblock In Van~Leeuwen, J., editor, {\em {Handbook of theoretical computer
  science (vol. A): algorithms and complexity}}. MIT Press.

\bibitem[\protect\citeauthoryear{von Neumann}{von Neumann}{1927}]{Neumann1927}
von Neumann, J. (1927).
\newblock Zur {H}ilbertschen {B}eweistheorie.
\newblock {\em Mathematische Zeitschrift\/}~{\bf 26}, 1--46.

\bibitem[\protect\citeauthoryear{Wang}{Wang}{1960}]{Wang1960}
Wang, H. (1960).
\newblock {Proving Theorems by Pattern Recognition I}.
\newblock {\em Communications of the ACM\/}~{\bf 3\/}(4), 220--234.

\bibitem[\protect\citeauthoryear{Wiedijk}{Wiedijk}{2006}]{Wiedijk2006}
Wiedijk, F. (2006).
\newblock {\em The seventeen provers of the world}, Volume 3600 of {\em Lecture
  Notes in Artificial Intelligence}.
\newblock Springer.

\bibitem[\protect\citeauthoryear{Williams}{Williams}{1998}]{Williams1998}
Williams, H. (1998).
\newblock {\em {\'E}douard Lucas and primality testing}.
\newblock Wiley.

\bibitem[\protect\citeauthoryear{Winker}{Winker}{1992}]{Winker1992}
Winker, S. (1992).
\newblock Absorption and idempotency criteria for a problem in near-boolean
  algebras.
\newblock {\em Journal of Algebra\/}~{\bf 153\/}(2), 414--423.

\end{thebibliography}

}

\end{document}